\documentclass[12pt]{article}

\usepackage{amssymb, amsmath, amsfonts, amsthm, mathrsfs, latexsym, epsfig,  psfrag, graphicx}
\usepackage[all]{xy}
\xyoption{v2} \xyoption{knot}\xyoption{arc} \xyoption{poly}
\xyoption{curve}

\usepackage{amsmath,amsthm,amsfonts,amssymb,amscd}
\def\qed{{\unskip\nobreak\hfil\penalty50
\hskip2em\hbox{}\nobreak\hfil$\square$
\parfillskip=0pt \finalhyphendemerits=0\par}\medskip}

\def\prf{\trivlist \item[\hskip \labelsep{\bf Proof\ }]}

\def\l{\langle}
\def\ro{\rho}

\def\ra{\rangle}
\def\Ad{{\mathrm {Ad}}}
\def\Aut{{\mathrm {Aut}}}
\def\Normal{{\mathrm {N}}}
\def\Auto{{\mathrm {Aut}}}
\def\Core{{\mathrm {Core}}}
\def\C{{\Bbb {C}}}
\def\End{{\mathrm {End}}}
\def\Exp{{\mathrm{Exp}}}
\def\Hom{{\mathrm {Hom}}}

\def\e{\varepsilon}

\def\la{\lambda}

\def\om{\omega}

\def\Ad{{\mathrm {Ad}}}
\def\Aut{{\mathrm {Aut}}}
\def\End{{\mathrm {End}}}
\def\Hom{{\mathrm {Hom}}}

\def\e{\varepsilon}

\def\la{\lambda}

\def\ro{{\rho}}

\newtheorem{theorem}{Theorem}[section]
\newtheorem{lemma}[theorem]{Lemma}
\newtheorem{conjecture}[theorem]{Conjecture}
\newtheorem{corollary}[theorem]{Corollary}
\newtheorem{definition}[theorem]{Definition}

\newtheorem{proposition}[theorem]{Proposition}
\newtheorem{remark}[theorem]{Remark}

\newtheorem{question}[theorem]{Question}
\newtheorem{notation}[theorem]{Notation}

\def\Hom{{\mathrm{Hom}}}

\def\res{\!\restriction\!}
\def\A{{\cal A}}

\def\C{{\mathbb C}}

\renewcommand{\qed}{\ \hfill $\blacksquare$}

\newcommand{\bdef}{\begin{definition}}
\newcommand{\blem}{\begin{lemma}}
\newcommand{\bprop}{\begin{proposition}}
\newcommand{\bthm}{\begin{theorem}}
\newcommand{\bcor}{\begin{corollary}}
\newcommand{\bconj}{\begin{conjecture}}
\newcommand{\ede}{\end{definition}}
\newcommand{\elem}{\end{lemma}}
\newcommand{\eprop}{\end{proposition}}
\newcommand{\ethm}{\end{theorem}}
\newcommand{\ecor}{\end{corollary}}
\newcommand{\econj}{\end{conjecture}}
\newcommand{\brem}{\begin{remark}}
\newcommand{\erem}{\end{remark}}

\newcommand{\ba}{\begin{array}}
\newcommand{\ea}{\end{array}}
\newcommand{\bea}{\begin{eqnarray}}

\title{\huge On intermediate subfactors of Goodman-de la Harpe-Jones subfactors\\}
\author{
{\sc Feng Xu}\footnote{Supported in part by NSF.}\\
Department of Mathematics\\
University of California at Riverside\\
Riverside, CA 92521\\
E-mail: {\tt xufeng@math.ucr.edu}}
\begin{document}
\date{}
\maketitle

\begin{abstract}
In this paper we present a conjecture on intermediate subfactors
which is a generalization of Wall's conjecture from the theory of
finite groups. Motivated by this conjecture, we determine all
intermediate subfactors of  Goodman-Harpe-Jones subfactors, and as a
result we verify that  Goodman-Harpe-Jones subfactors verify our
conjecture. Our result also gives a negative answer to a question
motivated by a conjecture of Aschbacher-Guralnick.

\end{abstract}

\newpage

\section{Introduction}

Let $M$ be a factor represented on a Hilbert space and $N$ a
subfactor of $M$ which is irreducible, i.e.,$N'\cap M= \C$. Let $K$
be an intermediate von Neumann subalgebra for the inclusion
$N\subset M.$ Note that $K'\cap K\subset N'\cap M = \C,$ $K$ is
automatically a factor. Hence the set of all intermediate subfactors
for $N\subset M$ forms a lattice under two natural operations
$\wedge$ and $\vee$ defined by:
\[
K_1\wedge K_2= K_1\cap K_2, K_1\vee K_2= (K_1\cup K_2)''.
\]
The commutant map $K\rightarrow K'$ maps an intermediate subfactor
$N\subset K\subset M$ to $M'\subset K'\subset N'.$ This map
exchanges the two natural operations defined above.

Let $M\subset M_1$ be the Jones basic construction of $N\subset M.$
Then  $M\subset M_1$ is canonically anti-isomorphic to $M'\subset N'$,
and the lattice of intermediate subfactors for $N\subset M$  is
related to the lattice of intermediate subfactors for $M\subset M_1$
by the commutant map defined as above.

 Let $G_1$ be a group and $G_2$ be a subgroup of $G_1$. An
interval sublattice $[G_1/G_2]$ is the lattice formed by all
intermediate subgroups $K, G_2\subseteq K\subseteq G_1.$

By cross product construction and Galois correspondence,  every
interval sublattice of finite groups can be realized as intermediate
subfactor lattice of finite index. Hence the study of intermediate
subfactor lattice of finite index is a natural generalization of the
study of interval sublattice of finite groups. The study of
intermediate subfactors has been very active in recent years (cf.
\cite{BJ},\cite{GJ}, \cite{JXu},\cite{Jh2}, \cite{ILP},
\cite{Longo4}, \cite{Wat} and \cite{ SW} for only a partial list).
There are a number of old problems about interval sublattice of
finite groups. It is therefore a natural programme to investigate if
these old problems have any generalizations to subfactor setting.
The hope is that maybe subfactor theory can provide new perspective
on these old problems.\par In \cite{Xl} we consider the problem
whether the very simple lattice $M_n$ consisting of a largest, a
smallest and $n$ pairwise incomparable elements can be realized as
subfactor lattice. We showed in \cite{Xl} all $M_{2n}$ are realized
as the lattice of intermediate subfactors of a pair of hyperfinite
type $III_1$ factors with finite depth. Since it is conjectured that
infinitely many $M_{2n}$ can not be realized as interval sublattices
of finite groups (cf. \cite{Balu} and \cite{Palfy}), our result
shows that if one is looking for obstructions for realizing finite
lattice as lattice of intermediate subfactors with finite index,
then the obstruction is very different from what one may find in
finite group theory.\par In 1961 G. E. Wall conjectured that the
number of maximal subgroups of a finite group $G$ is less than
$|G|$, the order of $G$ (cf. \cite{wall61}). In the same paper he
proved his conjecture when $G$ is solvable. See \cite{lie} for more
recent result on Wall's conjecture.

Wall's conjecture can be naturally generalized to a conjecture about
maximal elements in the lattice of intermediate subfactors. More
precisely, since $M$ is the maximal element under inclusion, what we
mean by maximal elements are those subfactors $K\neq M,N $ with the
property that if $K_1$ is an intermediate subfactor and $K\subset
K_1,$ then $K_1=M$ or $K.$ Minimal elements are defined similarly
where $N,M$ are not considered as minimal elements. When $M$ is the
cross product of $N$ by a finite group $G$, the maximal elements
correspond to maximal subgroups of $G,$ and the order of $G$ is the
dimension of second higher relative commutant. Hence a natural
generalization of Wall's conjecture is the following:
\begin{conjecture}\label{wall}
Let $N\subset M$ be an irreducible subfactor with finite index. Then
the number of maximal  intermediate subfactors is less than
dimension of $N'\cap M_1$ (the dimension of second higher relative
commutant of $N\subset M$).
\end{conjecture}
We note that since maximal intermediate subfactors in $N\subset M$
correspond to minimal intermediate subfactors in $M\subset M_1,$ and
the dimension of second higher relative commutant remains the same,
the conjecture  is equivalent to a similar conjecture as above with
maximal replaced by minimal.\par If we take $N$ and $M$ to be cross
products of a factor $P$ by $H$ and $G$ with $H$ a subgroup of $G$,
then the above conjecture gives a generalization of Wall's
conjecture which we call relative version of Wall's conjecture. The
relative version of Wall's conjecture states that the number of
maximal subgroups of $G$ strictly containing a subgroup $H$ is less
than $|G|/|H|.$

In the appendix we give a ``subfactor friendly" proof of relative
version of Wall's conjecture when $G$ is solvable. We also discuss a
question which is naturally motivated by a  conjecture  of
Aschbacher-Guralnick. This question also partially motivates our
work in this paper. A negative answer to this question is presented
in \S\ref{neg}.\par

When subfactors do not come from groups, with a few exceptions such as
\cite{GJ} and \cite{Xl}, very little is known about their maximal
intermediate subfactors. To test conjecture \ref{wall}, it is
therefore desirable to determine lattices of intermediate subfactors
for more examples of subfactors not coming from groups. As shown in
\cite{Xl}, a rich source of such subfactors come from conformal
field theories, and the techniques developed in \cite{Xl} allow one
to determine intermediate subfactors in many cases. In \cite{GJ}, as
part of effort to classify subfactors with no extra structure, the
intermediate subfactor lattice of a Goodman-Harpe-Jones (GHJ)
subfactor (cf. \cite{GHJ})  is determined. Since the dual GHJ
subfactors are closely related to conformal field theories based on
Loop group $LSU(2),$ we can use the method of \cite{Xl} to determine
lattices of intermediate subfactors for GHJ subfactors (The idea is
already presented in \cite{Xl}). In this paper we carry out this
idea. Compare to \cite{Xl}, the main difference  is that we need to
determine structure of a larger ring, but such rings have been
determined in \cite{Xb},\cite{BE3}, \cite{BEK2}. Combing these
results  we determine intermediate subfactor lattices for all GHJ
subfactors. One interesting consequence of our work is that the
intermediate subfactors of (dual) GHJ subfactors are again (dual)
GHJ subfactors. Also as a result we do not find any counter examples
to our conjecture \ref{wall}, and we  give a negative answer to
question \ref{gag}. We also find several surprising intermediate
subfactors which are not visible at first sight (cf. figures
\ref{d32}, \ref{d52}, \ref{e88}).

Since P. Grossman has proved that all non-commuting quadrilateral of
subfactors with bottom two subfactors of type $A$ come form GHJ
subfactor of type $D$ (cf. \cite{Gro}), figures \ref{dodd},
\ref{deven} also determined all intermediate subfactor lattices  of
non-commuting quadrilateral of subfactors with bottom two subfactors
of type $A$.

\par Besides what is already described
above, this paper is organized as follows: \S2 is a preliminary
section on sectors, representations of intermediate subfactors by a
pair of sectors, sectors from conformal nets, inductions,
Jones-Wassermann subfactors and a description of GHJ subfactors. The
simple idea of representations of intermediate subfactors by a pair
of sectors will prove to be crucial in later classifications. In \S3
we first explain the basic idea in \cite{Xl} to determine
intermediate subfactors by fusion, and  we carry out this idea for
GHJ subfactors of type $A$, $D$ and $E$ respectively. In \S4 we
apply the results of \S3 to determine the lattice relations of
intermediate subfactors, and these lattices are listed.

We'd like to thank Professors M. Aschbacher, P. Grossman and R.
Guralnick  for useful discussions, and especially Prof. V. F. R.
Jones for his interest and useful comments on conjecture \ref{wall}
which inspired the results of the appendix.
\section{Preliminaries}

For the convenience of the reader we collect here some basic notions
that appear in this paper. This is only a guideline and the reader
should look at the references such as preliminary sections of
\cite{KLX} for a more complete treatment.

\subsection{Sectors}
Let $M$ be a properly infinite factor and  $\text{\rm End}(M)$ the
semigroup of
 unit preserving endomorphisms of $M$.  In this paper $M$ will always
be the unique hyperfinite $III_1$ factors. Let $\text{\rm Sect}(M)$
denote the quotient of $\text{\rm End}(M)$ modulo unitary
equivalence in $M$. We  denote by $[\rho]$ the image of $\rho \in
\text{\rm End}(M)$ in  $\text{\rm Sect}(M)$.\par
 It follows from
\cite{Longo1} and \cite{Longo2} that $\text{\rm Sect}(M)$, with $M$
a properly infinite  von Neumann algebra, is endowed with a natural
involution $\theta \rightarrow \bar \theta $  ; moreover, $\text{\rm
Sect}(M)$ is
 a semiring.\par
Let $\rho \in \text{\rm End}(M)$ be a normal faithful conditional
expectation $\epsilon: M\rightarrow \rho(M)$.  We define a number
$d_\epsilon$ (possibly $\infty$) by:
$$
d_\epsilon^{-2} :=\text{\rm Max} \{ \lambda \in [0, +\infty)|
\epsilon (m_+) \geq \lambda m_+, \forall m_+ \in M_+ \}$$ (cf.
[PP]).\par
 We define
$$
d = \text{\rm Min}_\epsilon \{ d_\epsilon |  d_\epsilon < \infty \}.
$$   $d$ is called the statistical dimension of  $\rho$ and $d^2$ is called the Jones index of $\rho$.
It is clear from the definition that  the statistical dimension  of
$\rho$ depends only on the unitary equivalence classes  of  $\rho$.
The properties of the statistical dimension can be found in
\cite{Longo1},\cite{Longo2} and \cite{Longo3}.\par Denote by
$\text{\rm Sect}_0(M)$ those elements of $\text{\rm Sect}(M)$ with
finite statistical dimensions. For $\lambda $, $\mu \in \text{\rm
Sect}_0(M)$, let $\text{\rm Hom}(\lambda , \mu )$ denote the space
of intertwiners from $\lambda $ to $\mu $, i.e. $a\in \text{\rm
Hom}(\lambda , \mu )$ iff $a \lambda (x) = \mu (x) a $ for any $x
\in M$. $\text{\rm Hom}(\lambda , \mu )$  is a finite dimensional
vector space and we use $\langle \lambda , \mu \rangle$ to denote
the dimension of this space. $\langle  \lambda , \mu \rangle$
depends only on $[\lambda ]$ and $[\mu ]$. Moreover we have $\langle
\nu \lambda , \mu \rangle = \langle \lambda , \bar \nu \mu \rangle
$, $\langle \nu \lambda , \mu \rangle = \langle \nu , \mu \bar
\lambda \rangle $ which follows from Frobenius duality (See
\cite{Longo2} ).  We will also use the following notation: if $\mu $
is a subsector of $\lambda $, we will write as $\mu \prec \lambda $
or $\lambda \succ \mu $.  A sector is said to be irreducible if it
has only one subsector. \par For any $\rho\in \text{\rm End}(M)$
with finite index, there is a unique standard minimal inverse
$\phi_\rho: M\rightarrow M$ which satisfies
\[
\phi_\rho(\rho(m)m'\rho(m''))= m\phi_\rho(m')m'', m,m',m''\in M.
\]
$\phi_\ro$ is completely positive.  If $t\in \Hom(\rho_1,\rho_2)$
then we have
\begin{equation}\label{phiprop}
d_{\rho_1} \phi_{\rho_1} (mt)= d_{\rho_2} \phi_{\rho_2} (tm), m\in M
\end{equation}
\subsection{Representation of intermediate subfactors by a pair of
sectors} Let $M$ be an AFD type $III_1$ factor and $\rho\in
\End(M).$ Let $K$ be a factor such that $\rho(M)\subset K\subset M.$ Since $K$ is also AFD, one can
choose $\rho_1\in \End(M)$ with $\rho_1(M)=K.$ Then we have
$\rho=\rho_1\rho_2$ with $\rho_2=\rho_1^{-1}\rho\in {\rm End}(M).$
Conversely if $\rho=\rho_1\rho_2$ with $\rho_1,\rho_2 \in {\rm
End}(M),$ then $\rho(M)\subset \rho_1(M)\subset M.$ The
following lemma follows directly from definitions:
\begin{lemma}\label{pair}
Suppose that $\rho_1\rho_2= \sigma_1\sigma_2, \rho_i,\sigma_i\in
{\rm End}(M), i=1,2$ and $\rho_1(M)\subset \sigma_1(M)$. Then set
$\sigma=\sigma_1^{-1}\rho_1\in {\rm End}(M)$ we have
$\rho_1=\sigma_1\sigma, \sigma_2=\sigma\rho_2.$ Conversely if there
is $\sigma\in {\rm End}(M)$ such that $\rho_1=\sigma_1\sigma,
\sigma_2=\sigma\rho_2,$ then $\rho_1\rho_2= \sigma_1\sigma_2,$ and
$\rho_1(M)\subset \sigma_1(M).$ In addition $\sigma$ is an
automorphism iff $\rho_1(M)=\sigma_1(M).$
\end{lemma}
By Lemma \ref{pair} we can represent the intermediate subfactor of
$\rho$ by pairs $\rho_1,\rho_2$ with $\rho_1\rho_2=\rho,$ and if
$\sigma_1,\sigma_2$ represent the same intermediate subfactor iff
there is an automorphism $\sigma$ of $M$ such that
$\rho_1=\sigma_1\sigma, \sigma\rho_2=\sigma_2.$ The next lemma shows
that we can replace the pair $\rho_1,\rho_2$ by $[\rho_1],[\rho_2]$
when $\rho$ is irreducible:
\begin{lemma}\label{pairsect}
Suppose that $\rho=\rho_1\rho_2= \sigma_1\sigma_2,
\rho_i,\sigma_i\in {\rm End}(M)$, $[\rho_i]=[\sigma_i], i=1,2.$ Then
$\rho_1(M)=\sigma_1(M)$.
\end{lemma}
\prf By assumption we have unitaries $U_i\in M$ such that $\rho_i=
\Ad_{U_i}\sigma_i, i=1,2.$ Since  $\rho=\rho_1\rho_2=
\sigma_1\sigma_2,$ we have $\rho= \Ad_{U_1\sigma_1(U_2)}\rho.$ Since
$\rho$ is irreducible,  $U_1\sigma_1(U_2)$ must be a scalar multiple
of identity, and it follows that $\rho_1(M)= \Ad_{U_1}\sigma_1(M)
=\Ad_{\sigma_1(U_2^*)}\sigma_1(M)=
\sigma_1(\Ad_{U_2^*}(M))=\sigma_1(M) .$ \qed

In view of Lemma \ref{pair} and Lemma \ref{pairsect}, we introduce
the following notation:
\begin{definition}\label{pairdef}
We say that two pairs of sectors $[\rho_1], [\rho_2]$ and
$[\sigma_1],[\sigma_2]$ are equivalent if there is an automorphism
$\sigma$ of $M$ such that $[\rho_1]=[\sigma_1\sigma],
[\sigma\rho_2]=[\sigma_2].$ We denote the equivalence class of such
pair $[\rho_1], [\rho_2]$ by $[[\rho_1], [\rho_2]].$ when no
confusion arises we will write $[[\rho_1], [\rho_2]]$ simply as
$[\rho_1, \rho_2].$
\end{definition}
The following follows from our definition, Lemma \ref{pair} and
Lemma \ref{pairsect}:
\begin{corollary}\label{paircor}
Let $\rho\in {\rm End}(M)$ be irreducible. Then the set of
intermediate subfactors between $M$ and $\rho(M)$ can be
representanted naturally by $[[\rho_1], [\rho_2]]$ such that
$\rho_1\rho_2=\rho,$ and the intermediate subfactor is $\rho_1(M).$
In the following we will denote the intermediate subfactor
$\rho_1(M)$ simply by  $[[\rho_1], [\rho_2]]$. Then $[[\rho_1],
[\rho_2]]\subset [[\sigma_1], [\sigma_2]]$ as intermediate
subfactors between $M$ and $\rho(M)$ iff there is $\sigma\in {\rm
End}(M)$ such that $[\rho_1]=[\sigma_1\sigma],
[\sigma_2]=[\sigma\rho_2].$
\end{corollary}

\subsection{Symmetries of a subfactor}
\begin{definition}\label{autosub}
Let $\rho\in {\rm End}(M).$ Define $\Aut(\rho):=\{\sigma\in \Aut(M)|
\sigma\rho=\rho \}.$
\end{definition}

\blem\label{autosub1} If $\rho$ is irreducible and has finite index,
then $\Aut(\rho)$ is a finite group. There is a one to one
correspondence between  $\sigma\in\Aut(\rho)$ and sector $[\sigma]$
of index one which appears in $[\rho\bar\rho],$ and $\sigma$ is
constructed from $[\sigma]$ in the following way:  if $[\sigma]$ is
a sector of index one which appears in $[\rho\bar\rho],$ then we can
find a unique representative $\sigma\in \Aut(M)$ in the sector of
$[\sigma]$ such that $\sigma\rho=\rho.$

\elem \prf Let $\sigma\in \Aut(\rho).$ Since $\sigma\rho=\rho,$ by
Frobenius duality we have $\sigma\in \rho\bar\rho,$ and $[\sigma]$
is a sector of index one which appears in $[\rho\bar\rho]$ with
(necessarily) multiplicity one since $\rho$ is irreducible. One the
other hand if $[\sigma]$ is a sector of index one which appears in
$[\rho\bar\rho],$ then we can find a representative $\sigma$ in the
sector of $[\sigma]$ such that $\sigma\rho=\rho.$ Since $\rho$ is
irreducible, this  representative $\sigma$ in the sector of
$[\sigma]$ must be unique. Since $\rho$ has finite index, there are
only finitely many subsectors of $[\rho\bar\rho],$ and so
$\Aut(\rho)$ is finite. \qed

\begin{definition}\label{normalsub}
Let $\bar\rho\in {\rm End}(M).$ Denote by $U(M)$ (resp.
$U(\bar\rho(M))$ ) the group of unitaries in $M$ (resp.
$\bar\rho(M)$). Let $U_1(M)$ be the   subgroup of $U(M)$ which
consists of unitaries in $M$ whose conjugate action on $M$ preserve
$\bar\rho(M)$ as a set. Note that $U(\bar\rho(M))$ is a normal
subgroup of $U_1(M)$. Define
$\Normal(\bar\rho):=\frac{U_1(M)}{U(\bar\rho(M))}.$
\end{definition}

\blem\label{normalsub1} If $\rho$ is irreducible and has finite
index, then $\Normal(\bar\rho)$ is a finite group isomorphic to
$\Auto(\rho).$ In fact there is a one to one correspondence between
$U_\sigma\in\Normal (\bar\rho)$ and sector $[\sigma]$ of index one
which appears in $[\rho\bar\rho],$ and $U_\sigma$ is constructed
from $[\sigma]$ in the following way: if $[\sigma]$ is a sector of
index one which appears in $[\rho\bar\rho],$ then we can find a
unique representative $U_\sigma\in \Normal(\bar\rho)$  such that
$\Ad_{U_\sigma}\bar\rho\sigma=\bar\rho.$

\elem \prf Let $U$ be a unitary representative of an element in
$\Normal (\bar\rho).$ Since $\Ad_U\bar\rho=\bar\rho\sigma$ for some
automorphism $\sigma$ of $M$, by Frobenius duality we have
$\sigma\in \rho\bar\rho,$ and $[\sigma]$ is a sector of index one
which appears in $[\rho\bar\rho]$ with (necessarily) multiplicity
one since $\rho$ is irreducible. One the other hand if $[\sigma]$ is
a sector of index one which appears in $[\rho\bar\rho],$ then we can
find a unitary $U_\sigma$   such that
$\Ad_{U_\sigma}\bar\rho=\bar\rho\sigma.$ Since $\rho$ is
irreducible, the element with representative $U_\sigma$ in the group
$\Normal(\bar\rho)$  must be unique. Since $\rho$ has finite index,
there are only finitely many subsectors of $[\rho\bar\rho],$ and so
$\Normal (\bar\rho)$ is finite and it is isomorphic to $\Aut(\rho)$
by Lemma \ref{autosub1}.
 \qed

\subsection{Sectors from conformal nets and their representations}
We refer the reader to \S3 of \cite{KLX} for definitions of
conformal nets and their representations. Suppose a conformal net
$\A$ and a representation $\lambda$ is given. Fix an open interval I
of the circle and Let $M:=\A(I)$ be a fixed type $III_1$ factor.
Then $\lambda$ give rises to an endomorphism still denoted by
$\lambda$ of $M$. We will recall some of the results of \cite{R2}
and introduce notations. \par Suppose $\{[\lambda] \}$ is a finite
set  of all equivalence classes of irreducible, covariant,
finite-index representations of an irreducible local conformal net
$\A$. We will use $\Delta_\A$ or simply $\Delta$ to denote all
finite index representations of net $\A$ and will use the same
notation $\Delta_\A$ to denote the corresponding sectors of $M.$

We will denote the conjugate of $[\lambda]$ by $[{\bar \lambda}]$
and identity sector (corresponding to the vacuum representation) by
$[0]$ if no confusion arises, and let $N_{\lambda\mu}^\nu = \langle
[\lambda][\mu], [\nu]\rangle $. Here $\langle \mu,\nu\rangle$
denotes the dimension of the space of intertwiners from $\mu$ to
$\nu$ (denoted by $\text {\rm Hom}(\mu,\nu)$).  We will denote by
$\{T_e\}$ a basis of isometries in $\text {\rm
Hom}(\nu,\lambda\mu)$. The univalence of $\lambda$ and the
statistical dimension of (cf. \S2  of \cite{GL1}) will be denoted by
$\omega_{\lambda}$ and $d{(\lambda)}$ (or $d_{\lambda})$)
respectively. The unitary braiding operator $\epsilon (\mu,
\lambda)$ (cf. \cite{GL1} ) verifies the following
\begin{proposition}\label{bfe}
\text{\rm (1)}  Yang-Baxter-Equation (YBE)
$$
\e(\mu , \gamma ) \mu (\e(\lambda , \gamma )) \e(\lambda , \mu ) =
\gamma (\e(\lambda , \mu )) \e(\lambda , \gamma )\lambda (\e(\mu ,
\gamma ))\, .
$$

\text{\rm (2)}  Braiding-Fusion-Equation (BFE)

For any $w\in \text{\rm Hom} (\mu \gamma , \delta )$
\begin{align*} \e(\lambda , \delta )\lambda (w) = w\mu (\e(\lambda , \gamma ))
\e(\lambda , \mu ) \\
\e(\delta , \lambda )w = \lambda (w) \e(\mu , \lambda ) \mu
(\e(\gamma , \lambda )) \, \\
\e( \delta,\lambda )^* \lambda (w) = w\mu (\e( \gamma,\lambda )^*)
\e(\mu,\lambda )^* \\
\e(\lambda, \delta )^* \lambda (w) = w\mu (\e( \gamma,\lambda )^*)
\e(\lambda,\mu )^*
\end{align*}
\end{proposition}

\blem\label{monodromy} If $\lambda,\mu$ are irreducible, and
$t_\nu\in \Hom(\nu, \lambda\mu)$ is an isometry, then $t_\nu
\e(\mu,\lambda)\e(\lambda,\mu) t_\nu^* =
\frac{\om_\nu}{\om_\lambda\om_\mu}.$ \elem

By Prop. \ref{bfe}, it follows that if $t_i\in \Hom (\mu_i,\lambda)$
is an isometry, then
\[
\e(\mu,\mu_i)\e(\mu_i,\mu)=t_i^*\e(\mu,\la)\e(\lambda,\mu) t_i
\]
We shall always identify the center of $M$ with $\C.$  Then we have
the following

\blem\label{escalar} If
\[\e(\mu,\la)\e(\lambda,\mu)\in \C,\]
then \[\e(\mu,\mu_i)\e(\mu_i,\mu)\in \C, \forall \mu_i\prec \la. \]
\elem

Let $\phi_\lambda$ be the unique minimal left inverse of $\lambda$,
define:
\begin{equation}\label{ymatrix}
Y_{\lambda\mu}:= d(\lambda)  d(\mu) \phi_\mu (\epsilon (\mu, \lambda)^*
\epsilon (\lambda, \mu)^*),
\end{equation}
where $\epsilon (\mu, \lambda)$ is the unitary braiding operator
 (cf. \cite{GL1} ). \par
We list two properties of $Y_{\lambda \mu}$ (cf. (5.13), (5.14) of
\cite{R2}):
\begin{lemma}\label{yprop}
\begin{equation*}
Y_{\lambda\mu} = Y_{\mu\lambda}  = Y_{\lambda\bar \mu}^* =
Y_{\bar \lambda \bar \mu}.
\end{equation*}
\begin{equation*}
Y_{\lambda\mu}  = \sum_k N_{\lambda\mu}^\nu \frac{\omega_\lambda\omega_\mu}
{\omega_\nu} d(\nu) .
\end{equation*}
\end{lemma}
We note that one may take the second equation in the above lemma as the
definition of $Y_{\lambda\mu}$.\par
Define
$a := \sum_i d_{\rho_i}^2 \omega_{\rho_i}^{-1}$.
If the matrix $(Y_{\mu\nu})$ is invertible,
by Proposition on P.351 of \cite{R2} $a$ satisfies
$|a|^2 = \sum_\lambda d(\lambda)^2$.
\begin{definition}\label{c0'}
Let $a= |a| \exp(-2\pi i \frac{c_0}{8})$ where  $c_0\in {\mathbb R}$ and $c_0$
is well defined ${\rm mod} \ 8\mathbb Z$.
\end{definition}
Define matrices
\begin{equation}\label{Smatrix}
S:= |a|^{-1} Y, T:=  C {\rm Diag}(\omega_{\lambda})
\end{equation}
where \[C:= \label{dims} \exp(-2\pi i \frac{c_0}{24}).\]
Then these matrices satisfy (cf. \cite{R2}):
\begin{lemma}\label{Sprop}
\begin{align*}
SS^{\dag} & = TT^{\dag} ={\rm id},  \\
STS &= T^{-1}ST^{-1},  \\
S^2 & =\hat{C},\\
 T\hat{C} & =\hat{C}T,
\end{align*}

where $\hat{C}_{\lambda\mu} = \delta_{\lambda\bar \mu}$
is the conjugation matrix.
\end{lemma}
Moreover
\begin{equation}\label{Verlinde}
N_{\lambda\mu}^\nu = \sum_\delta \frac{S_{\lambda\delta} S_{\mu\delta}
S_{\nu\delta}^*}{S_{1\delta}}. \
\end{equation}
is known as Verlinde formula.
The commutative algebra  generated by $\lambda$'s with structure
constants $N_{\lambda\mu}^\nu$ is called {\bf fusion algebra} of
$\A$. If $Y$ is invertible,
it follows from Lemma \ref{Sprop}, (\ref{Verlinde})  that any nontrivial
irreducible representation
of the fusion algebra is of the form
$\lambda\rightarrow \frac{S_{\lambda\mu}}{S_{1\mu}}$ for some
$\mu$. \par
\subsection{Induced endomorphisms}\label{inductionsection}
Suppose that $\ro\in \End(M)$ has the property that
$\gamma=\rho\bar\rho \in \Delta_\A$. By \S2.7 of \cite{LR}, we can
find two isometries $v_1\in \Hom(\gamma,\gamma^2), w_1\in \Hom
(1,\gamma)$\footnote{We use $v_1,w_1$ instead of $v,w$ here since
$v, w$ are used to denote sectors in Section \ref{typea}.} such that
$\bar\rho(M)$ and $v_1$ generate $M$ and
\begin{align*}
v_1^* w_1 & = v_1^* \gamma(w_1) = d_\rho^{-1} \\
v_1v_1 & = \gamma(v_1) v_1  \\
\end{align*}
By Thm. 4.9 of \cite{LR}, we shall say that $\rho$ is {\it local }
if
\begin{align}\label{local}
v_1^* w_1 & = v_1^* \gamma(w_1) = d_\rho^{-1} \\
v_1v_1 & = \gamma(v_1) v_1  \\
\bar\rho(\epsilon(\gamma,\gamma)) v_1 & = v_1
\end{align}
Note that if $\rho$ is local, then
\begin{equation}\label{local=1}
\om_\mu =1, \forall \mu\prec \ro\bar\ro \end{equation}

For each (not necessarily irreducible)  $\lambda\in \Delta_\A,$ let
$\e{(\lambda,\gamma)}
 :   \lambda\gamma \rightarrow
\gamma\lambda$ (resp.  $\tilde\e{(\lambda,\gamma)}$), be the
positive (resp. negative) braiding operator as defined in Section
1.4 of \cite{Xb}. Denote by $\lambda_\e \in$ End$(M)$ which is
defined by
\begin{align*}
\lambda_\e(x) :&= ad(\e(\lambda,\gamma))\lambda(x)=
\e(\lambda,\gamma) \lambda(x) \e(\lambda,\gamma)^* \\
\lambda_{\tilde\e}(x) :&= ad(\tilde\e(\lambda,\gamma))\lambda(x)=
\tilde\e(\lambda,\gamma)^* \lambda(x) \tilde\e(\lambda,\gamma)^*,
\forall x\in M.
\end{align*}
 By (1) of Theorem 3.1 of \cite{Xb},  $\lambda_{\e} \rho (M)
\subset \rho(M), \lambda_{\tilde\e}  \rho (M) \subset \rho(M),$
hence the following definition makes sense\footnote{We have changed
the notations $a_\la, \tilde a_\la$ of \cite{Xb} to $\tilde a_\la,
a_\la$ of this paper to make some of the formulas such as equation
(\ref{aa2}) simpler.}. \bdef\label{ala} If $\la\in \Delta_\A$ define
two elements of $\End(M)$ by
$$
a^\ro_\lambda(m):= \rho^{-1} (\lambda_\e  \rho (m)) , \ \tilde
a^\ro_\lambda(m):= \rho^{-1} (\lambda_{\tilde\e}  \rho (m)), \forall
m\in M.
$$
$a_\lambda^\ro$ (resp. $\tilde a_\lambda^\ro$) will be referred to
as positive (resp. negative) induction of $\la$ with respect to
$\ro.$

\ede
\begin{remark}
For simplicity we will use $a_\la, \tilde a_\la$ to denote
$a_\lambda^\ro, \tilde a_\lambda^\ro$ when it is clear that
inductions are with respect to the same $\ro.$
\end{remark}

The endomorphisms $a_\lambda$ are called braided endomorphisms in
\cite{Xb}  due to its braiding properties (cf.  (2) of Corollary 3.4
in \cite{Xb}), and enjoy an interesting set of properties (cf.
Section 3 of \cite{Xb}). Though \cite{Xb} focus on the local case
which was clearly the most interesting case in terms of producing
subfactors, as observed in
 \cite{BE1}, \cite{BE2}, \cite{BE3}, \cite{BE4}
 that many of the arguments in \cite{Xb} can be
generalized. These properties are also studied in a slightly
different context in \cite{BE1}, \cite{BE2}, \cite{BE3}. In these
papers,  the induction is between $M$ and a subfactor $N$ of $M$
,while the induction above is on the same algebra. A dictionary
between our notations here and these papers has been set up in
\cite{Xj} which simply use an isomorphism between $N$ and $M$. Here
one has a choice to use this isomorphism to translate all
endomorphisms of $N$ to endomorphims of $M$, or equivalently all
endomorphims of $M$ to endomorphims of $N$. In \cite{Xj} the later
choice is made (Hence $M$ in \cite{Xj} will be our $N$ below). Here
we make the first choice which makes the dictionary slightly
simpler. Our dictionary here is equivalent to that of \cite{Xj}. Set
$N=\bar\rho(M)$. In the following the notations from \cite{BE1} will
be given a subscript BE. The formulas are :
\begin{align}\label{aa1}
\rho\res N=i_{BE}, \ & \bar\ro\ro\res N=\bar i_{BE}i_{BE}, \\
\gamma = \bar\rho^{-1} \theta_{BE} \bar\rho, \ &
\bar\ro \ro= \gamma_{BE}, \\
\lambda = \bar\ro^{-1} \lambda_{BE} \bar\ro, \ \ & \e(\lambda, \mu)
= \bar\rho (\e^{+}(\lambda_{BE}, \mu_{BE})) \\
& \tilde\e(\lambda, \mu) = \bar\rho (\e^{-}(\lambda_{BE}, \mu_{BE}))
\end{align}
The dictionary between $a_\lambda \in End (M)$ in definition
\ref{ala} and $\alpha_\lambda^{-}$ as in Definition 3.3, Definition
3.5 of \cite{BE1}  are given by:
\begin{equation}\label{aa2}
a_\lambda  =\alpha^+_{\lambda_{BE}}, \tilde a_\lambda
=\alpha^{-}_{\lambda_{BE}}
\end{equation}
The above formulas will be referred to as our {\it dictionary}
between the notations of \cite{Xb} and that of \cite{BE1}. The proof
is the same as that of \cite{Xj}. Using this dictionary one can
easily translate results of \cite{Xb} into the settings of
\cite{BE1},  \cite{BE2}, \cite{BE3}, \cite{BE4},\cite{BEK1},
\cite{BEK2}and vice versa. First we summarize a few properties from
\cite{Xb} which will be used in this paper: (cf. Th. 3.1 , Co. 3.2
and Th. 3.3 of \cite{Xb} ): \bprop\label{xua} (1). The maps
$[\lambda] \rightarrow [a_\lambda], [\lambda]\rightarrow [\tilde
a_\lambda]$ are ring homomorphisms;\par (2) $ a_\lambda \bar \rho
=\tilde a_\lambda \bar \rho =\bar \rho \lambda$;\par (3) When
$\ro\bar\ro$ is local, $ \l a_\la,  a_\mu \ra = \l \tilde a_\la,
\tilde a_\mu \ra = \l a_\la\bar \ro, a_\mu\bar\ro \ra =\l \tilde
a_\la\bar \ro, \tilde a_\mu\bar\ro \ra ;$\par (4) (3) remains valid
if $a_\la, a_\mu$ (resp. $\tilde a_\la, \tilde a_\mu$) are replaced
by their subsectors.

\eprop

\bdef\label{hro} $H_\ro$ is a finite dimensional vector space over
$\C$ with orthonormal basis consisting of irreducible sectors of
$[\lambda\rho],\forall \lambda\in \Delta_\A.$ \ede
\par

$[\lambda]$ acts linearly on $H_\ro$ by $[\lambda] [a]= \sum_b \l
\lambda a,b\ra [b]$ where $[b]$ are elements in the basis of
$H_\ro.$ \footnote{By abuse of notation, in this paper we use
$\sum_b$ to denote the sum over the basis $[b]$ in $H_\ro$.}  By
abuse of notation, we use $[\lambda]$ to denote the corresponding
matrix relative to the basis of $H_\ro.$ By definition these
matrices are normal and commuting, so they can be simultaneously
diagonalized. Recall the irreducible representations of the fusion
algebra of $\A$ are given by
$$
\lambda \rightarrow \frac{S_{\lambda \mu}}{S_{1\mu}}.
$$
\bdef\label{orthogonalphig} Assume $\l \la a,b\ra=
\sum_{{\mu,i}\in(\text{\rm Exp})} \frac{S_{\lambda \mu}}{S_{0\mu}}
\cdot \phi_a^{(\mu,i)} \phi_b^{(\mu,i)^*}$ where $ \phi_a^{(\mu,i)}$
are normalized orthogonal eigenvectors of $[\lambda]$  with
eigenvalue $\frac{S_{\lambda \mu}}{S_{0\mu}},$ $\Exp$ is a set of
$\mu,i$'s
 and $i$ is an
index indicating the multiplicity of  $\mu$, and is called the set
of exponents of $H_\rho.$  Recall if a representation is denoted by
$0$, it will always be the vacuum representation. \ede

The following lemma is elementary:

\blem\label{ephi} (1):
\[\sum_b d_b^2=\frac{1}{S_{00}^2}\]
where the sum is over the basis of $H_\ro$. The vacuum appears once
in $\Exp$ and
\[
\phi_a^{(1)}= S_{00} d_a; \]\par (2) \[\sum_i
\frac{\phi_a^{(\la,i)}{{\phi_b^{(\la,i)}}^*}}{S_{0\la}^2}= \sum_\nu
\l \bar\nu a,b\ra \frac{S_{\nu\la}}{S_{0\la}}\] where if $\la$ does
not appear in $\Exp$ then the righthand side is zero.

\elem

\prf Ad (1): By definition we have
\[[a\bar\rho]= \sum_{\la}\l a\bar\rho,\la\ra[\la]=\sum_{\la}\l a,\la\rho\ra[\la]
\] where in the second $=$ we have used Frobenius reciprocity. Hence
\[d_a d_{\bar\ro}= \sum_\lambda \l a\bar\rho,\la\ra d_\la\] and we
obtain
\[
\sum_\la d_\la^2 = \sum_{\la,a} \l a\bar\ro,\la\ra d_\la d_a/d_\ro =
\sum_a d_a^2
\]
(2) follows from definition and orthogality of $S$ matrix. \qed

\par

In \cite{BE3} and \cite{Xb}, commutativity among subsectors of
$a_\la,\tilde a_\mu, \la,\mu\in \Delta$ were studied.  We record
these results in the following for later use: \blem\label{acommut}
(1) Let $[b]$ (resp. $[b']$) be any subsector of $a_\la$ (resp.
$\tilde a_\la$). Then
\[
[a_\mu b]=[ba_\mu], [\tilde a_\mu b']=[b'\tilde a_\mu]\forall \mu,
[bb']=[bb'];
\]\par
(2) Let $[b]$ be a subsector of $a_\mu\tilde a_\la,$ then $[a_\nu
b]=[b a_\nu],[\tilde a_\nu b]=[b \tilde a_\nu], \forall \nu;$\par
(3) If $\sigma\prec \la \rho,$ the $[\mu \sigma]=[\sigma
a_\mu]=[\sigma \tilde{a_\mu}].$
\elem


\subsection {Jones-Wassermann subfactors from representation of Loop
groups}\label{typea} Let $G= SU(n)$. We denote $LG$ the group of
smooth maps $f: S^1 \mapsto G$ under pointwise multiplication. The
diffeomorphism group of the circle $\text{\rm Diff} S^1 $ is
naturally a subgroup of $\text{\rm Aut}(LG)$ with the action given
by reparametrization. In particular the group of rotations
$\text{\rm Rot}S^1 \simeq U(1)$ acts on $LG$. We will be interested
in the projective unitary representation $\pi : LG \rightarrow U(H)$
that are both irreducible and have positive energy. This means that
$\pi $ should extend to $LG\rtimes \text{\rm Rot}\ S^1$ so that
$H=\oplus _{n\geq 0} H(n)$, where the $H(n)$ are the eigenspace for
the action of $\text{\rm Rot}S^1$, i.e., $r_\theta \xi = \exp(i n
\theta)$ for $\theta \in H(n)$ and $\text{\rm dim}\ H(n) < \infty $
with $H(0) \neq 0$. It follows from \cite{PS} that for fixed level
$k$ which is a positive integer, there are only finite number of
such irreducible representations indexed by the finite set
$$
 P_{++}^{k}
= \bigg \{ \lambda \in P \mid \lambda = \sum _{i=1, \cdots , n-1}
\lambda _i \Lambda _i , \lambda _i \geq 0\, , \sum _{i=1, \cdots ,
n-1} \lambda _i \leq k \bigg \}
$$
where $P$ is the weight lattice of $SU(n)$ and $\Lambda _i$ are the
fundamental weights. We will write $\la=(\la_1,...,\la_{n-1}),
\la_0= k-\sum_{1\leq i\leq n-1} \la_i$ and refer to
$\la_0,...,\la_{n-1}$ as components of $\la.$

We will use $\Lambda_0$ or simply $1$  to denote the trivial
representation of $SU(n)$. For $\lambda , \mu , \nu \in P_{++}^{k}$,
define $N_{\lambda \mu}^\nu  = \sum _{\delta \in P_{++}^{k}
}S_\lambda ^{(\delta)} S_\mu ^{(\delta)} S_\nu
^{(\delta*)}/S_{\Lambda_0}^{(\delta})$ where $S_\lambda ^{(\delta)}$
is given by the Kac-Peterson formula:
$$
S_\lambda ^{(\delta)} = c \sum _{w\in S_n} \varepsilon _w \exp
(iw(\delta) \cdot \lambda 2 \pi /n)
$$
where $\varepsilon _w = \text{\rm det}(w)$ and $c$ is a
normalization constant fixed by the requirement that
$S_\mu^{(\delta)}$ is an orthonormal system. It is shown in
\cite{Kac2} P. 288 that $N_{\lambda \mu}^\nu $ are non-negative
integers. Moreover, define $ Gr(C_k)$ to be the ring whose basis are
elements of $ P_{++}^{k}$ with structure constants $N_{\lambda
\mu}^\nu $.
  The natural involution $*$ on $ P_{++}^{k}$ is
defined by $\lambda \mapsto \lambda ^* =$ the conjugate of $\lambda
$ as representation of $SU(n)$.\par

We shall also denote $S_{\Lambda _0}^{(\Lambda)}$ by $S_1^{(\Lambda
)}$. Define $d_\lambda = \frac {S_1^{(\lambda )}}{S_1^{(\Lambda
_0)}}$. We shall call $(S_\nu ^{(\delta )})$ the $S$-matrix of
$LSU(n)$ at level $k$. \par
 The
following  result is proved in \cite{Wass} (See Corollary 1 of
Chapter V in \cite{Wass}).

\bthm\label{wass}  Each $\lambda \in  P_{++}^{(k)}$ has finite index
with index value $d_\lambda ^2$.  The fusion ring generated by all
$\lambda \in P_{++}^{(k)}$ is isomorphic to $ Gr(C_k)$. \ethm

\begin{remark}\label{jxx}
The subfactors in the above theorem are called Jones-Wassermann
subfactors after the authors who first studied them (cf.
\cite{Jh1},\cite{Wass}).
\end{remark}
We will concentrate on $n=2$ case in this paper. Fix level $k\geq
1,$ the representations are labeled by half integers $i$ with $0\leq
i\leq k/2.$ For example $0$ will label the vacuum representation,
and $1/2$ will label the vector representation. The statistical
dimension of $1/2$ is $2\cos(\frac{2\pi}{k+2}).$ The fusion rules
are given by
$$
[i][j]=\sum_{|i-j|\leq l\leq i+j, i+j+l\leq k,i+j+l\in {\Bbb Z}}[l].
$$
The $S$ matrix for $n=2$ is given by $S_{ij}=\sqrt{\frac{2}{k+2}}
\sin(\frac{\pi(2i+1)(2j+1)}{k+2})$.

\subsection{Dual Goodman-Harpe-Jones subfactors from Conformal Field Theory}
The dual GHJ subfactors are obtained from irreducible sectors of
$[i][\rho]$ with $\rho\bar\rho\in \Delta_\A$ (cf.
\cite{Xb},\cite{BEK2}, Appendix of \cite{BE3}).  The induction in
the following will be with respect to $\rho.$
\begin{definition}\label{fusiongraph}
The fusion graph of $[1/2]$ on $H_\rho$ is the graph whose vertices
are irreducible sectors in $H_\rho$ and two vertices $a,b$ are
connected by an edge if $\langle 1/2 a,b\rangle=1.$ We will say that
$H_\rho$ is type $A, D, E$ if the fusion graph is $A, D, E$
respectively. By an end point of the fusion graph we mean a vertex
on the graph which is connected to only one other point on the
graph.
\end{definition}
Since the norm of the fusion graph of $[1/2]$ is less than 2, the
graph must be $A-D-E$ (cf. \cite{GHJ}). When $\rho=0,$  the GHJ
subfactors are Jones subfactors. The fusion graph of $1/2$ is given
by:

\begin{figure}[h]
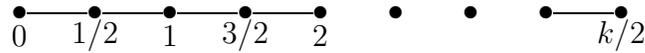

\[
\xy (0,0)*{\bullet}; (10,0)*{\bullet}**\dir{-};
(20,0)*{\bullet}**\dir{-}; (30,0)*{\bullet}**\dir{-};
(40,0)*{\bullet}**\dir{-}; (50,0)*{\bullet}; (60,0)*{\bullet};
(70,0)*{\bullet};(80,0)*{\bullet}**\dir{-};
 (0,-3)*{0}; (10,-3)*{{1/2}};
(20,-3)*{{1}};(30,-3)*{ {3/2}};(40,-3)*{ {2}}; (80,-3)*{k/2}
\endxy
\]
\caption{Fusion graph of $1/2$, type $A_{k+1}$}\label{a}
\end{figure}

When $k$ is even, and $[\rho\bar\rho]=[0]+[k/2].$ The fusion graph
of $[1/2]$ on $H_\rho$ is type $D$ graph. In this case $[i][\rho]$
is irreducible iff $i\neq k/4$. If $4|k$ then $[a_{k/4}]=[b]+[b'],$
and if $[\bar\rho\rho]=[0]+[g],[g^2]=[0],$ then $[g{b}g]=[{b'}].$ If
$k$ is not divisible by $4,$ $[k/4][\rho]=[b]+[b']$ where $b,b'$ are
irreducible with same index, and $[\bar\rho\rho]=[0]+[a_{k/2}]$ and
$[k/2][b]=[b']$ (cf. \S5.2 of \cite{BEK2}).
\begin{notation}
For convenience we shall use $\rho_0$ to denote a fixed endomorphism
with $[\rho_0\bar\rho_0]=[0]+[k/2].$
\end{notation}
\par
\begin{figure}[h]
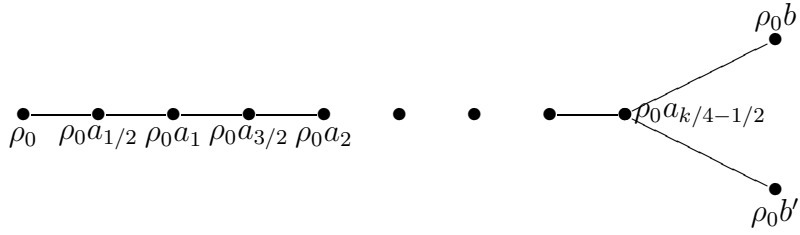

\[
\xy (0,0)*{\bullet}; (10,0)*{\bullet}**\dir{-};
(20,0)*{\bullet}**\dir{-}; (30,0)*{\bullet}**\dir{-};
(40,0)*{\bullet}**\dir{-}; (50,0)*{\bullet}; (60,0)*{\bullet};
(70,0)*{\bullet};(80,0)*{\bullet}**\dir{-};(100,10)*{\bullet}**\dir{-};
(80,0)*{\bullet};(100,-10)*{\bullet}**\dir{-};
 (0,-3)*{\rho_0}; (10,-3)*{\rho_0 a_{1/2}};
(20,-3)*{\rho_0 a_{1}};(30,-3)*{\rho_0 a_{3/2}};(40,-3)*{\rho_0
a_{2}}; (90,0)*{\rho_0 a_{k/4-1/2}};(100,13)*{\rho_0
b};(100,-13)*{\rho_0 b'}
\endxy
\]
\caption{Fusion graph of $1/2$, type $D_{\frac{k}{2}+2},
4|k$}\label{deven1}
\end{figure}

\begin{figure}[h]
\[
\xy (0,0)*{\bullet}; (10,0)*{\bullet}**\dir{-};
(20,0)*{\bullet}**\dir{-}; (30,0)*{\bullet}**\dir{-};
(40,0)*{\bullet}**\dir{-}; (50,0)*{\bullet}; (60,0)*{\bullet};
(70,0)*{\bullet};(80,0)*{\bullet}**\dir{-};(100,10)*{\bullet}**\dir{-};
(80,0)*{\bullet};(100,-10)*{\bullet}**\dir{-};
 (0,-3)*{\rho_0}; (10,-3)*{\rho_0 a_{1/2}};
(20,-3)*{\rho_0 a_{1}};(30,-3)*{\rho_0 a_{3/2}};(40,-3)*{\rho_0
a_{2}}; (90,0)*{\rho_0 a_{k/4-1/2}};(100,13)*{ b};(100,-13)*{b'}
\endxy
\]
\caption{Fusion graph of $1/2$, type $D_{\frac{k}{2}+2}, 4\not |k
$}\label{dodd1}
\end{figure}

\begin{figure}[h]
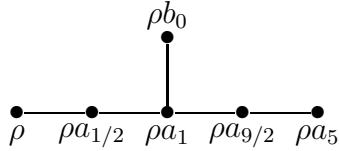

\[
\xy (0,0)*{\bullet}; (10,0)*{\bullet}**\dir{-};
(20,0)*{\bullet}**\dir{-}; (30,0)*{\bullet}**\dir{-};
(40,0)*{\bullet}**\dir{-}; (20,0)*{}; (20,10)*{\bullet}**\dir{-};
(0,-3)*{\rho}; (10,-3)*{\rho a_{1/2}}; (20,-3)*{\rho
a_{1}};(30,-3)*{\rho a_{9/2}};(40,-3)*{\rho a_{5}};
(20,13)*{\rho
b_0}
\endxy
\]
\caption{Fusion graph of $1/2$, type $E_6$}\label{e6}
\end{figure}

When the fusion graph of $[1/2]$ is $E_6,$ the graph is given by
figure \ref{e6} (cf. Page 392-393 of \cite{Xb}).

Note that  $k=10, [b_0]=[a_{3/2}]-[a_{9/2}],
[\rho\bar\rho]=[0]+[3].$

\begin{figure}[h]
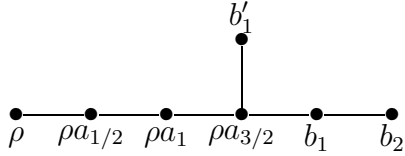

\[
\xy (0,0)*{\bullet}; (10,0)*{\bullet}**\dir{-};
(20,0)*{\bullet}**\dir{-}; (30,0)*{\bullet}**\dir{-};
(40,0)*{\bullet}**\dir{-};(50,0)*{\bullet}**\dir{-}; (30,0)*{};
(30,10)*{\bullet}**\dir{-}; (0,-3)*{\rho}; (10,-3)*{\rho a_{1/2}};
(20,-3)*{\rho a_{1}};(30,-3)*{\rho a_{3/2}};(40,-3)*{ b_1};
(30,13)*{ b_1'};(50,-3)*{ b_2}
\endxy
\]
\caption{Fusion graph of $1/2$, type $E_7$}\label{e7}
\end{figure}

When the fusion graph of $[1/2]$ is $E_7,$ $k=16,
[\rho\bar\rho]=[0]+[4]+[8],$ and the graph is given by figure
\ref{e7}(cf. Figure 42 of \cite{BEK2}).
Note that
$$
[b_2]=[5/2\rho]-[3/2\rho], [b_1]=[1/2b_2], [b_1']=[5/2\rho]-[b_1]
$$

\begin{figure}[h]
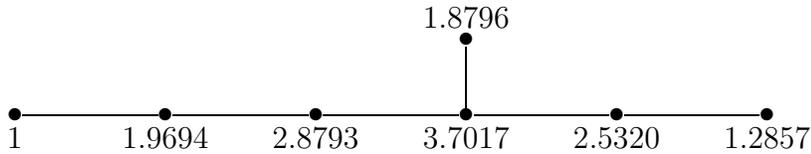

\[
\xy (0,0)*{\bullet}; (20,0)*{\bullet}**\dir{-};
(40,0)*{\bullet}**\dir{-}; (60,0)*{\bullet}**\dir{-};
(80,0)*{\bullet}**\dir{-};(100,0)*{\bullet}**\dir{-}; (60,0)*{};
(60,10)*{\bullet}**\dir{-}; (0,-3)*{1}; (20,-3)*{1.9694};
(40,-3)*{2.8793};(60,-3)*{3.7017};(80,-3)*{2.5320}; (60,13)*{
1.8796};(100,-3)*{ 1.2857}
\endxy
\]
\caption{Normalized eigenvectors}\label{e7eig}
\end{figure}

The square root of indices of the corresponding subfactors divided
by $d_\rho$ is given by figure \ref{e7eig}.

\begin{figure}[h]
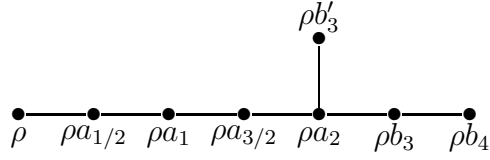

\[
\xy (0,0)*{\bullet}; (10,0)*{\bullet}**\dir{-};
(20,0)*{\bullet}**\dir{-}; (30,0)*{\bullet}**\dir{-};
(40,0)*{\bullet}**\dir{-};(50,0)*{\bullet}**\dir{-};(60,0)*{\bullet}**\dir{-};
(40,0)*{}; (40,10)*{\bullet}**\dir{-}; (0,-3)*{\rho}; (10,-3)*{\rho
a_{1/2}}; (20,-3)*{\rho a_{1}};(30,-3)*{\rho a_{3/2}}; (40,-3)*{
\rho a_2}; (50,-3)*{ \rho b_3}; (40,13)*{\rho b_3'};(60,-3)*{ \rho
b_4}
\endxy
\]
\caption{Fusion graph of $1/2$, type $E_8$}\label{e8}
\end{figure}

\begin{figure}[h]
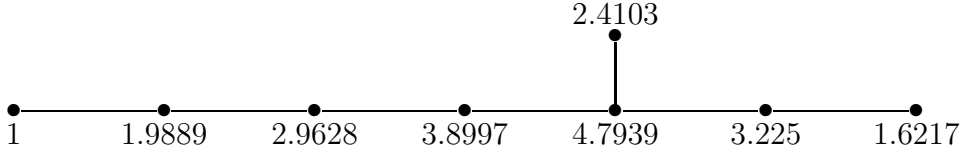

\[
\xy (0,0)*{\bullet}; (20,0)*{\bullet}**\dir{-};
(40,0)*{\bullet}**\dir{-}; (60,0)*{\bullet}**\dir{-};
(80,0)*{\bullet}**\dir{-};(100,0)*{\bullet}**\dir{-};(120,0)*{\bullet}**\dir{-};
(80,0)*{}; (80,10)*{\bullet}**\dir{-}; (0,-3)*{1}; (20,-3)*{1.9889};
(40,-3)*{2.9628};(60,-3)*{3.8997}; (80,-3)*{4.7939};
(100,-3)*{3.225}; (80,13)*{2.4103};(120,-3)*{1.6217}
\endxy
\]
\caption{Normalized eigenvectors}\label{e8eig}
\end{figure}

\begin{figure}[h]
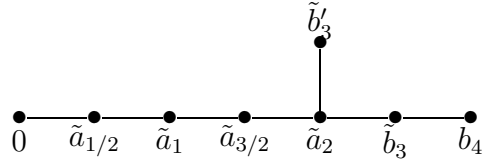

\[
\xy (0,0)*{\bullet}; (10,0)*{\bullet}**\dir{-};
(20,0)*{\bullet}**\dir{-}; (30,0)*{\bullet}**\dir{-};
(40,0)*{\bullet}**\dir{-};(50,0)*{\bullet}**\dir{-};(60,0)*{\bullet}**\dir{-};
(40,0)*{}; (40,10)*{\bullet}**\dir{-}; (0,-3)*{0}; (10,-3)*{\tilde
a_{1/2}}; (20,-3)*{\tilde a_{1}};(30,-3)*{\tilde a_{3/2}}; (40,-3)*{
\tilde a_2}; (50,-3)*{ \tilde b_3}; (40,13)*{\tilde b_3'};(60,-3)*{
 b_4}
\endxy
\]
\caption{Fusion graph of $\tilde a_{1/2}$}\label{til}
\end{figure}

When the fusion graph of $[1/2]$ is $E_8,$ $k=28,
[\rho\bar\rho]=[0]+[5]+[9]+[14]$ (cf. \cite{Xb}),  and the graph is
given by figure \ref{e8}.
The square root of indices of the corresponding subfactors divided
by $d_\rho$ is given by figure \ref{e8eig}.
The fusion graph of $\tilde a_{1/2}$ is given by figure \ref{til}.
Note that we have
$$
[a_{5/2}]=[b_3']+[b_3], [b_3]=[a_{1/2} b_4], [\tilde
a_{5/2}]=[\tilde b_3']+[\tilde b_3], [\tilde b_3]=[\tilde a_{1/2}
b_4].
$$
\section{Classification of intermediate subfactors of dual GHJ subfactors}
Let $\sigma$ be a dual GHJ subfactor. We assume that $k\geq 5$ in
this section. Recall that a subfactor is maximal if there is no
nontrivial intermediate subfactor. To list nontrivial intermediate
subfactors of $\sigma$, according to Cor. \ref{paircor}, we need to
determine all pairs $[\sigma_1,\sigma_2]$ such that
$[\sigma_1\sigma_2]=[\sigma]$ with $1<d_{\sigma_1}<d_{\sigma}.$
Since $\sigma_1\bar\sigma_1\prec \sigma\bar\sigma\in \Delta,$ by
Definition \ref{ala} we can apply induction with respect to
$\sigma_1.$

Our basic strategy, as explained in the end of \S2 of \cite{Xl}, is
to consider the fusion graph of $[1/2]$ on $H_{\sigma_1}.$ By
$A-D-E$ classification, this graph must be one of $A-D-E$ graphs.
The following is useful: \blem\label{end} Let $\sigma$ be a dual GHJ
subfactor such that $\sigma=\sigma_1\sigma_2.$ Then either $[1/2
\sigma_1]$ or both $[a_{1/2}\sigma_2]$ and $[\tilde
a_{1/2}\sigma_2]$ are both irreducible where the inductions are with
respect to $\sigma_1$. In terms of fusion graphs of $[1/2]$ on
$\sigma_1,$ $[a_{1/2}]$ on $\sigma_2$ and $[\tilde a_{1/2}]$ on
$\sigma_2,$ this means that either $\sigma_1$ is an end point or
$\sigma_2$ is an end point on the  fusion graphs of $[a_{1/2}]$ on
$\sigma_2$ and $[\tilde a_{1/2}]$ on $\sigma_2.$

\elem \prf Since $\sigma$ is irreducible, we have $ 1=\langle
\sigma_1\sigma_2, \sigma_1\sigma_2\rangle =\langle
\bar\sigma_1\sigma_1, \sigma_2\bar\sigma_2\rangle $ On the other
hand
$$
\langle 1/2\sigma_1, 1/2\sigma_1\rangle =\langle \sigma_1 a_{1/2},
\sigma_1 a_{1/2}\rangle =\langle \bar\sigma_1\sigma_1,
a_{1/2}^2\rangle =1+\langle \bar\sigma_1\sigma_1, a_{1}\rangle
$$
$$
\langle 1/2\sigma_1, 1/2\sigma_1\rangle =\langle \sigma_1 \tilde
a_{1/2}, \sigma_1 \tilde a_{1/2}\rangle =\langle
\bar\sigma_1\sigma_1, \tilde a_{1/2}^2\rangle =1+\langle
\bar\sigma_1\sigma_1, \tilde a_{1}\rangle
$$
Similarly
$$
\langle a_{1/2}\sigma_2, a_{1/2}\sigma_2\rangle = \langle
\sigma_2\bar\sigma_2, a_{1/2}^2\rangle =1+\langle
\sigma_2\bar\sigma_2, a_{1}\rangle
$$
$$
\langle \tilde a_{1/2}\sigma_2, \tilde a_{1/2}\sigma_2\rangle =
\langle \sigma_2\bar\sigma_2, \tilde a_{1/2}^2\rangle =1+\langle
\sigma_2\bar\sigma_2, \tilde a_{1}\rangle
$$
Since $a_1, \tilde a_1$ are irreducible by Lemma 2.33 of \cite{Xl},
the lemma follows. \qed

By Lemma \ref{end}, either $\sigma_1$ or $\sigma_2$ is an end point
on a graph of type $A-D-E,$ which has two or three end points. This
together with the known values of indices greatly reduce possible
choices of $\sigma_1,\sigma_2,$ and in the next few sections we will
determine all such choices.
\subsection{Type A}
By Example 5.24 of \cite{Xl}, in this case $i$ is maximal if $i\neq
k/4.$ When $k=2,3$ all GHJ subfactors are maximal, and when $k=4$ a
direct computation gives the lattice of intermediate subfactors
represented by sector  $[1]$ as in figure \ref{k4}. Assume now
$k\geq 5.$
\begin{lemma}\label{speca}
Assume that $[k/4]=[\sigma_1][\sigma_2].$ Then if $0\leq j\leq k/2$
and $j$ is an integer, then $j\in \Exp$ of $H_{\sigma_1}$ as defined
in definition \ref{orthogonalphig}.
\end{lemma}
\prf Apply (2) of Lemma \ref{ephi} to $a=\bar\sigma_2, b=\sigma_1,$
we have
\[\sum_i
\frac{\phi_a^{(j,i)}{{\phi_b^{(j,i)}}^*}}{S_{0j}^2}= \sum_\nu \l
\bar\nu a,b\ra \frac{S_{\nu j}}{S_{0j}} =\frac{S_{k/4j}}{S_{0j}}
\]
Since $\frac{S_{k/4j}}{S_{0j}}$ up to a nonzero constant is
$\sin(\frac{(2j+1)\pi}{2})\neq 0$ when $0\leq j\leq k/2$ and $j$ is
an integer, it follows that  $j\in \Exp$ of $H_{\sigma_1}$ as
defined in definition \ref{orthogonalphig}. \qed

When $k$ is even and $i=k/4,$ there are two cases to consider: when
$4|k,$ by Lemma \ref{speca} if $j$ is an integer then $j\in \Exp,$
and it follows that $H_{\sigma_1}$ is either type $A$ or $D$. When
$H_{\sigma_1}$ is type $A$ then by fusion rule up to right
multiplication by an automorphism $[\sigma_1]=[\sigma]$ or
$d_{\sigma_1}=1;$ When $H_{\sigma_1}$ is type $D,$ by left
multiplication by an automorphism on $\sigma_1$ if necessary we may
assume that $\rho\in H_{\sigma_1},$ with $[\rho\bar\rho]=[0]+[k/2].$
By Lemma \ref{end} either $\sigma_1$ or $\sigma_2$ has to be an end.
It follows that $[\sigma_1,\sigma_2]=[\rho{b},\bar\rho]$ or
$[\sigma_1,\sigma_2]=[\rho,{b}\bar\rho]. $

\par

When $k$ is not divisible by $4,$ similar argument shows that
$[\sigma_1,\sigma_2]=[{b},\bar\rho]$ or
$[\sigma_1,\sigma_2]=[\rho,\bar b].$

We summarize the result in the following
\begin{theorem}\label{casea}
Suppose that $k\geq 4.$ If $[i]\neq [k/4],$ the corresponding
subfactor is maximal; If $k=4,$ the intermediate subfactor of $[1]$
is given by $[\rho {b},\bar\rho]$  $[\rho\bar\rho]=[0]+[1],$ $\rho
{b}\prec 2\rho;$ If $4|k,k\geq 8$, the intermediate subfactors of
$[k/4]$ is given by $[\rho {b},\bar\rho]$ and $[\rho,{b}\bar\rho]$
where $[\rho\bar\rho]=[0]+[k/2],$ $\rho{b}\prec k/2\rho;$ If $k$ is
even but not divisible by 4, the intermediate subfactors of $[k/4]$
is given by $[b,\bar\rho]$ and $[\rho,\bar{b}]$ where
$[\rho\bar\rho]=[0]+[k/2],$ ${b}\prec k/2\rho.$
\end{theorem}

\subsection{Type $E_6$}
In this section we assume that $\sigma$ is a dual GHJ subfactor
appearing in $j\rho$ with $[\rho\bar\rho]=[0]+[3].$ Let
$\sigma=\sigma_1\sigma_2.$\par If $H_{\sigma_1}$ is type $A$, then
multiply $\sigma_1$ on the right by an automorphism if necessary, we
may assume that $\sigma_1\in \Delta,$ hence $\sigma_1\in
\Delta\rho.$ By Lemma \ref{end} either $1/2\sigma_1$ is irreducible
or $a_{1/2}\sigma_2$ is irreducible. Then from figure \ref{e6} it is
clear that the possible pairs of $[\sigma_1,\sigma_2]$ are given by
$$
[1/2,\rho],[1,\rho],[9/2,\rho],[4,\rho], [1/2,\rho b]
$$
If  $H_{\sigma_1}$ is type $D$, we may assume that $\sigma_1\in
\Delta \rho_0,$ and it follows that $\sigma_2\in \bar\rho_0
\Delta\rho.$ We note that $[5\rho b_0]=[\rho b_0],[5\rho a_1]=[\rho
a_1].$ By Lemma \ref{autosub}  $\rho b_0(M) $ (resp. $\rho a_1(M)$)
is contained in a subfactor of index $2$ which is the fixed point
subalgebra of $M$ under an automorphism determined by $[5].$ Hence
there are sectors $x,y$ such that $[\rho_0 x]=[\rho b_0], [\rho_0
y]=[\rho a_1].$ Since $k=10,$ we have
$[\bar\rho_0\rho_0]=[0]+[a^{\rho_0}_5].$ From $ 1=\langle \rho_0
x,\rho_0 x\rangle= 1 +\langle  a^{\rho_0}_5 x, x\rangle$ we conclude
that $[x]\neq [ a^{\rho_0}_5 x].$ Similarly
 $[y]\neq [ a^{\rho_0}_5 y].$
Let us determine all irreducible sectors of $\bar\rho_0 \Delta\rho.$
Note that $a_{1/2}^{\rho_0}$ acts on $\sigma_2\in \bar\rho_0
\Delta\rho,$ and the corresponding fusion graph of
$a_{1/2}^{\rho_0}$ is an $A-D-E$ graph. Since $\langle
\bar\rho_0\rho, \bar\rho_0\rho\rangle= \langle \rho_0\bar\rho_0,
\rho\bar\rho\rangle =1,$  $ \langle a_{1/2}^{\rho_0}\bar\rho_0\rho,
a_{1/2}^{\rho_0}\bar\rho_0\rho\rangle= \langle
([0]+[1])\rho_0\bar\rho_0, \rho\bar\rho\rangle =1,$ and
$[a_1^{\rho_0}\bar\rho_0\rho] = [\bar\rho_0\rho y]
=[y]+[a^{\rho_0}_5 y],$ it follows that the fusion graph of the
action of $a_{1/2}^{\rho_0}$ is given by (The vertices are labeled
by all irreducible sectors of $\bar\rho_0 \Delta\rho$) figure
\ref{e6d}.
\begin{figure}[h]
\[
\xy (0,0)*{\bullet}; (10,0)*{\bullet}**\dir{-};
(20,0)*{\bullet}**\dir{-}; (30,0)*{\bullet}**\dir{-};
(40,0)*{\bullet}**\dir{-}; (20,0)*{}; (20,10)*{\bullet}**\dir{-};
(0,-3)*{x}; (10,-3)*{y};
(20,-6)*{a_{1/2}^{\rho_0}\bar\rho_0\rho};(30,-3)*{ a^{\rho_0}_5
y};(40,-3)*{a^{\rho_0}_5x}; (20,13)*{\bar\rho_0\rho}
\endxy
\]
\caption{Fusion graph of $a_{1/2}^{\rho_0}$}\label{e6d}
\end{figure}
By Lemma \ref{end} and figure \ref{e6d}, it follows that the
following is a list of possible intermediate subfactors:
$$
[\rho_0,x], [\rho_0,y], [1/2\rho_0,x]
$$
Note that from $[\rho_0 x]=[\rho b_0]$ we have
$$
[\rho \bar x]\prec [\rho \bar b_0 \bar\rho \rho_0] = [\rho]
([a_{3/2}]-[a_{9/2}]) [\bar\rho \rho_0]
=([3/2]-[9/2])([0]+[3])[\rho_0] =2[3/2\rho_0]
$$
Similarly
$$
[\rho \bar x a^{\rho_0}_5]\prec [\rho \bar b_0 \bar\rho \rho_0] =
[\rho] ([a_{3/2}]-[a_{9/2}])[ \bar\rho \rho_0]
=([3/2]-[9/2])([0]+[3])[\rho_0] =2[3/2\rho_0]
$$
By comparing indices   we have proved the following
\begin{equation}\label{3/2}
[\rho \bar{x}]=[\rho \bar{x} a^{\rho_0}_5]=[3/2\rho_0]
\end{equation}
When $H_{\sigma_1}$ is $E_6,$ then there is $\rho_2\in H_{\sigma_1}$
such that $[\rho_2\bar\rho_2]=[0]+[3].$ By the cohomology vanishing
result of remark 5.4 in \cite{KLV}, multiplying $\sigma_1$ on the
right by an automorphism if necessary, we can assume that
$[\rho_2]=[\rho],$ and so $\sigma_1\in \Delta \rho,$ and
$\sigma_2\in \bar \rho\Delta \rho.$ The set of irreducible sectors
of $\bar \rho\Delta \rho$ and fusion graphs of the action of
$a_{1/2}, \tilde a_{1/2}$ are given by \cite{Xb} and Figure 5 of
\cite{BE3}. By using Lemma \ref{end} it is straightforward to
determine the following list of possible intermediate subfactors:
$$
[\rho, a_{1/2}], [\rho, \tilde a_{1/2}], [\rho, a_{9/2}], [\rho,
\tilde a_{9/2}], [\rho, a_{1}], [\rho, \tilde a_{1}],[\rho, b],
[\rho b,  a_{1/2}], [\rho b, \tilde a_{1/2}], [1/2\rho, b].
$$

\begin{theorem}\label{casee6}
Assume that $\sigma$ is a dual GHJ subfactor appearing in $j\rho$
with $[\rho\bar\rho]=[0]+[3].$ Then the following is the list of all
intermediate subfactors that can occur:
$$
[1/2,\rho],[1,\rho],[9/2,\rho],[4,\rho], [1/2,\rho b], [\rho_0,x],
[\rho_0,y], [1/2\rho_0,x], [\rho, a_{1/2}], [\rho, \tilde a_{1/2}],
[\rho, a_{9/2}]
$$

$$[\rho, \tilde a_{9/2}], [\rho, a_{1}], [\rho,
\tilde a_{1}],[\rho, b], [\rho b,  a_{1/2}], [\rho b, \tilde
a_{1/2}], [1/2\rho, b].
$$
\end{theorem}

\subsection{$E_7$ case}
In this section we assume that $\sigma$ is a dual GHJ subfactor
appearing in $j\rho$ with $[\rho\bar\rho]=[0]+[4]+[8].$ Since
$[8\rho]=[\rho],$ by Lemma \ref{autosub} we may assume that
$\rho=\rho_o\rho_1,$ and similarly
$$
[b_1]=[\rho_0\hat{b_1}], [b_1']=[\rho_0\hat{b_1'}],
[b_2]=[\rho_0\hat{b_2}].
$$

Let $\sigma=\sigma_1\sigma_2.$ If $H_{\sigma_1}$ is type $A$, then
multiply $\sigma_1$ on the right  by an automorphism if necessary,
we may assume that $\sigma_1\in \Delta,$ hence $\sigma_1\in
\Delta\rho.$ By Lemma \ref{end} either $1/2\sigma_1$ is irreducible
or $a_{1/2}\sigma_2$ is irreducible. Then from figure \ref{e7} it is
clear that the possible pairs of $[\sigma_1,\sigma_2]$ are given by
$$
[1/2,\rho],[1,\rho],[3/2,\rho],[1/2,b_1'], [1/2, b_2], [1,b_2].
$$
If  $H_{\sigma_1}$ is type $D$, we may assume that $\sigma_1\in
\Delta \rho_0,$ and it follows that $\sigma_2\in \bar\rho_0
\Delta\rho.$ From $\rho=\rho_0\rho_1$ we get
$[\bar\rho_0\rho]=[\rho_1]+[g\rho_1]$ where
$[\bar\rho_0\rho_0]=[0]+[g], [g^2]=[0].$ It is then easy to
determine all irreducible sectors of $\bar\rho_0 \Delta \rho.$ There
are two $E_7$ graphs in figure \ref{e7d} encoding the action of
$a_{1/2}^{\rho_0}.$
\begin{figure}[h]
\[
\xy (0,0)*{\bullet}; (10,0)*{\bullet}**\dir{-};
(20,0)*{\bullet}**\dir{-}; (30,0)*{\bullet}**\dir{-};
(40,0)*{\bullet}**\dir{-};(50,0)*{\bullet}**\dir{-}; (30,0)*{};
(30,10)*{\bullet}**\dir{-}; (0,-3)*{\rho_1}; (10,-3)*{\rho_1
a_{1/2}}; (20,-3)*{\rho_1 a_{1}};(30,-3)*{\rho_1 a_{3/2}};(40,-3)*{
\hat{b_1}}; (30,13)*{ \hat{b_1'}};(50,-3)*{ \hat{b_2}}
\endxy
\]
\[
\xy (0,0)*{\bullet}; (20,0)*{\bullet}**\dir{-};
(40,0)*{\bullet}**\dir{-}; (60,0)*{\bullet}**\dir{-};
(80,0)*{\bullet}**\dir{-};(100,0)*{\bullet}**\dir{-}; (60,0)*{};
(60,10)*{\bullet}**\dir{-}; (0,-3)*{g\rho_1}; (20,-3)*{g\rho_1
a_{1/2}}; (40,-3)*{g\rho_1 a_{1}};(60,-3)*{g\rho_1
a_{3/2}};(80,-3)*{ g\hat{b_1}}; (60,13)*{ g\hat{b_1'}};(100,-3)*{
g\hat{b_2}}
\endxy
\]
\caption{Fusion graphs of $a_{1/2}^{\rho_0}$}\label{e7d}
\end{figure}

By Lemma \ref{end} and figure \ref{e7d} it follows that the
following is a list of possible intermediate subfactors:
$[\rho_0,w]$ where $w$ is one of the vertices in the first graph of
figure \ref{e7d}, $ [\rho_0 b, \hat{b_2}], [\rho_0 b, g\hat{b_2}], $
where $\rho_0b\prec [4\rho_0], [\rho_0 b \hat{b_2}]=[\rho_0 b
g\hat{b_2}]=[3/2\rho],$
$$
[1/2\rho_0,\rho_1],[1\rho_0,\rho_1],[3/2\rho_0,\rho_1],[1/2\rho_0,\hat{b_1'}],[1/2\rho_0,\hat{b_2}],[1\rho_0,\hat{b_2}]
$$

Note that from $[5/2\rho]=[3/2\rho]+[b_2]$ we have
$$
[5/2\rho\bar\rho_1]=[3/2\rho\bar\rho_1]+[b_2\bar\rho_1],
[5/2\rho\bar\rho_1g]=[3/2\rho\bar\rho_1g]+[b_2\bar\rho_1g]
$$
On the other hand from $[\rho]=[\rho_0\rho_1]$ and computing index
it is easy to derive that
$$
[\rho_1\bar\rho_1]+[g\rho_1\bar\rho_1 g]=2[0]+[a_4^{\rho_0}]
$$
or
$$
[\rho_1\bar\rho_1]+[g\rho_1\bar\rho_1 g]=2[0]+[\tilde a_4^{\rho_0}]
$$
Plug into the equations above and comparing index we find that
$$
[{b_2}\bar\rho_1]=[b_2\bar\rho_1g]=[5/2\rho_0]
$$
It follows that
$$
[\hat{b_2}\bar\rho_1]\prec
[\bar\rho_05/2\rho_0]=[a_{5/2}^{\rho_0}]+[ga_{5/2}^{\rho_0}]
$$
Hence $[\hat{b_2}\bar\rho_1]=[ga_{5/2}^{\rho_0}]$ or
$[\hat{b_2}\bar\rho_1]=[a_{5/2}^{\rho_0}].$ By taking conjugates we
have $[\hat{b_2}\bar\rho_1]=[\rho_1\bar{\hat{b_2}}].$ We have
therefore proved the following:
\begin{equation}\label{5/2}
[{b_2}\bar\rho_1]=[b_2\bar\rho_1g]=[5/2\rho_0]=[\rho\bar{\hat{b_2}}]=[\rho\bar{\hat{b_2}}g]
\end{equation}

When $H_{\sigma_1}$ is $E_7,$ then there is $\rho_2\in H_{\sigma_1}$
such that $[\rho_2\bar\rho_1]=[\rho\bar\rho]].$ By the cohomology
vanishing result of remark 5.4 in  \cite{KLV}, multiplying
$\sigma_1$ on the right by an automorphism if necessary, we can
assume that $[\rho_2]=[\rho],$ and so $\sigma_1\in \Delta \rho,$ and
$\sigma_2\in \bar \rho\Delta \rho.$ The set of irreducible sectors
of $\bar \rho\Delta \rho$ and fusion graphs of the action of
$a_{1/2}, \tilde a_{1/2}$ are given by  Figure 42 of \cite{BEK2}. By
using Lemma \ref{end} it is straightforward to determine the
following list of possible intermediate subfactors:
$$
[\rho, a_{1/2}], [\rho, \tilde a_{1/2}], [\rho, a_{3/2}], [\rho,
\tilde a_{3/2}], [\rho, a_{1}], [\rho, \tilde a_{1}],
[\rho,\tau],[b_1',a_{1/2}],[b_1',\tilde
a_{1/2}],[b_2,a_{1/2}],[b_2,\tilde a_{1/2}]
$$
$$
[b_2,a_{1}],[b_2,\tilde a_{1}], [b_2,\delta]
$$
where $[\tau]=[a_{1/2}]([\tilde a_{5/2}]-[\tilde a_{3/2}]),
[\delta]=[a_4]-[\tilde a_1]$ and
$$
[\rho\tau]=[b_1], [b_2\delta]=[3/2\rho].
$$

\begin{theorem}\label{casee7}
Assume that $\sigma$ is a dual GHJ subfactor appearing in $j\rho$
with $[\rho\bar\rho]=[0]+[4]+[8].$ Then the following is the list of
all intermediate subfactors that can occur:
$$
[1/2,\rho],[1,\rho],[3/2,\rho],[1/2,b_1'], [1/2, b_2], [1,b_2]
$$

$[\rho_0,w]$ where $w$ is one of the vertices in the first graph of
figure \ref{e7d},
$$
[\rho_0 b, \hat{b_2}], [\rho_0 b, g\hat{b_2}],
$$ where $\rho_0b\prec [4\rho_0], [\rho_0 b \hat{b_2}]=[\rho_0 b
g\hat{b_2}]=[3/2\rho],$
$$
[1/2\rho_0,\rho_1],[1\rho_0,\rho_1],[3/2\rho_0,\rho_1],[1/2\rho_0,\hat{b_1'}],[1/2\rho_0,\hat{b_2}],[1\rho_0,\hat{b_2}]
$$
$$
[\rho, a_{1/2}], [\rho, \tilde a_{1/2}], [\rho, a_{3/2}], [\rho,
\tilde a_{3/2}], [\rho, a_{1}], [\rho, \tilde a_{1}],
[\rho,\tau],[b_1',a_{1/2}],[b_1',\tilde
a_{1/2}],[b_2,a_{1/2}],[b_2,\tilde a_{1/2}]
$$
$$
[b_2,a_{1}],[b_2,\tilde a_{1}], [b_2,\delta]
$$
where $[\tau]=[a_{1/2}]([\tilde a_{5/2}]-[\tilde a_{3/2}]),
[\delta]=[a_4]-[\tilde a_1]$ and
$$
[\rho\tau]=[b_1], [b_2\delta]=[3/2\rho].
$$

\end{theorem}

\subsection{Type D}
In this section we assume that $\sigma$ is a dual GHJ subfactor
appearing in $j\rho$ with $[\rho\bar\rho]=[0]+[k/2].$ Let
$\sigma=\sigma_1\sigma_2.$
\begin{lemma}\label{specd}
Assume that $\sigma$ is a dual GHJ subfactor appearing in $j\rho$
with $[\rho\bar\rho]=[0]+[k/2].$ Let $\sigma=\sigma_1\sigma_2.$ If
$0\leq i\leq k/2$ and $i$ is an integer, and $\sin(\frac{\pi
(2i+1)(2j+1)}{k+2})\neq 0,$ then $i\in \Exp$ of $H_{\sigma_1}$ as
defined in definition \ref{orthogonalphig}.
\end{lemma}
\prf When $j\neq k/4,$ we have $[\sigma]=[j\rho].$ Apply (2) of
Lemma \ref{ephi} to $a=\rho\bar\sigma_2, b=\sigma_1,$ we have
\[\sum_l
\frac{\phi_a^{(i,l)}{{\phi_b^{(i,l)}}^*}}{S_{0i}^2}= \sum_\nu \l
\bar\nu a,b\ra \frac{S_{\nu i}}{S_{0i }} =\frac{S_{j i}}{S_{0i}}+
\frac{S_{(k/2-j) i}}{S_{0i}}
\]
When $j= k/4,$ we have $[\sigma\bar\rho]=[k/4].$ Apply (2) of Lemma
\ref{ephi} to $a=\rho\bar\sigma_2, b=\sigma_1,$ we have
\[\sum_l
\frac{\phi_a^{(i,l)}{{\phi_b^{(i,l)}}^*}}{S_{0i}^2}= \sum_\nu \l
\bar\nu a,b\ra \frac{S_{\nu i}}{S_{0i }} =\frac{S_{\frac{k}{2}
i}}{S_{0i}}
\]

Note that  $\frac{S_{ji}}{S_{0i}}+\frac{S_{(k/2-j) i}}{S_{0i}}=0 $
if $i$ is not an integer. When $i$ is an integer,  up to a nonzero
constant $\frac{S_{ji}}{S_{0i}}+\frac{S_{(k/2-j) i}}{S_{0i}} $ is
$\sin(\frac{\pi (2i+1)(2j+1)}{k+2}),$ it follows that if $0\leq
i\leq k/2,$ $i$ is an integer, and $\sin(\frac{\pi
(2i+1)(2j+1)}{k+2})\neq 0,$ then $i\in \Exp$ of $H_{\sigma_1}$ as
defined in definition \ref{orthogonalphig}. \qed

By Lemma \ref{specd}, if $i$ is an integer, $(2i+1)(2j+1)$ is not
divisible by $k+2,$ then $i\in \Exp$ of $H_{\sigma_1}.$ By
inspecting exponents on Page 18 of \cite{GHJ}, if $k\neq 10,16,$
$H_{\sigma_1}$ must be of type $A$ or $D.$ We will see in the
following that when $k=10,16$ $H_{\sigma_1}$ can be $E_6, E_7$
respectively.
\subsubsection{Local case: $4|k, k\neq 16$}

In this section we assume that $4|k, k\neq 16 $ and $\sigma$ is a
dual GHJ subfactor appearing in $j\rho$ with
$[\rho\bar\rho]=[0]+[k/2].$ By our assumption $H_{\sigma_1}$ must be
of type $A$ or $D.$ If $H_{\sigma_1}$ is type $A$, then multiply
$\sigma_1$ on the right by an automorphism if necessary, we may
assume that $\sigma_1\in \Delta,$ hence $\sigma_1\in \Delta\rho.$ By
Lemma \ref{end} either $1/2\sigma_1$ is irreducible or
$a_{1/2}\sigma_2$ is irreducible. It is clear that the possible
pairs of $[\sigma_1,\sigma_2]$ are given by
$$
[i,\rho],i\neq [k/4],[1/2,\rho b],[1/2,\rho b'],
$$
where $[a_{k/4}=[b]+[b'], [1/2\rho b]=[1/2\rho b']=[\rho
a_{k/4-1/2}].$\par If $H_{\sigma_1}$ is type $D$, we may assume that
$\sigma_1\in \Delta \rho,$ and it follows that $\sigma_2\in \bar\rho
\Delta\rho.$ It follows that the fusion graph of the action of
$a_{1/2}^{\rho_0}$ is given by two $D$ graphs, whose vertices are
labeled by irreducible components of $a_{i}, i\in \Delta$ and
$\tilde a_{i}, i\in \Delta$ respectively. By Lemma \ref{end}, it
follows that the following is a list of possible intermediate
subfactors: $ [\rho,w] $ where $w$ is an irreducible component of
$a_{i},i\in \Delta,$
$$
[\rho b, a_{1/2}], [\rho b,\tilde a_{1/2}], [1/2\rho,b],
[1/2\rho,b']
$$
\begin{theorem}\label{keven}
Assume that $4|k, k\neq 16 $ and $\sigma$ is a  dual GHJ subfactor
appearing in $j\rho$ with $[\rho\bar\rho]=[0]+[k/2].$ Then the
following is a list of possible intermediate subfactors:
$$
[i,\rho],i\neq [k/4],[1/2,\rho b],[1/2,\rho b'],
$$
where $[a_{k/4}=[b]+b'], [1/2\rho b]=[1/2\rho b']=[\rho
a_{k/4-1/2}]$
 $ [\rho,w] $ where $w$ is an irreducible component of
$a_{i},i\in \Delta,$
$$
[\rho b, a_{1/2}], [\rho b,\tilde a_{1/2}], [1/2\rho,b],
[1/2\rho,b']
$$
\end{theorem}
\subsubsection {$k=16$}
By Lemma \ref{specd}, $H_{\sigma_1}$ may be $E_7$ when $\sigma=
5/2\rho_0=\sigma_1\sigma_2.$ Suppose that $H_{\sigma_1}$ is $E_7.$
Then we can assume that $\sigma_1\in \Delta \rho_2$ with
$[\rho_2\bar\rho_2]=[0]+[4]+[8].$ It follows that $\sigma_2\in
\bar\rho_2 \Delta \rho_0.$ But the conjugate of $\bar\rho_2 \Delta
\rho_0$ has already been determined by Figure \ref{e7d}. By
computing index we conclude that either $d_{\sigma_1}=
d_{b_2},d_{\sigma_2}= d_{\rho_1}$ or $d_{\sigma_1}=
d_{\rho_2},d_{\sigma_2}= d_{\hat{b_2}}.$  By equation \ref{5/2} we
conclude that $[\sigma_1,\sigma_2]=[b_2,\bar\rho_1],
[\sigma_1,\sigma_2]=[b_2,\bar\rho_1 g],$
$[\sigma_1,\sigma_2]=[\rho_2,\bar{\hat {b_2}}] $, or
$[\sigma_1,\sigma_2]=[\rho_2,\bar{\hat {b_2}}g] $.
\begin{theorem}\label{16}
Assume that $k=16$ and $\sigma$ is a dual GHJ subfactor appearing in
$j\rho$ with $[\rho\bar\rho]=[0]+[8].$ Then the following is a list
of possible intermediate subfactors:
$$
[i,\rho],i\neq [k/4],[1/2,\rho b],[1/2,\rho b'],
$$
where $[a_{k/4}]=[b]+b'], [1/2\rho b]=[1/2\rho b']=[\rho
a_{k/4-1/2}]$
 $ [\rho,w] $ where $w$ is an irreducible component of
$a_{i},i\in \Delta,$
$$
[\rho b, a_{1/2}], [\rho b,\tilde a_{1/2}], [1/2\rho,b],
[1/2\rho,b']
$$
$$
[b_2,\bar\rho_1], [b_2,\bar\rho_1 g], [\rho_2,\bar{\hat {b_2}}],
[\rho_2,\bar{\hat {b_2}}g]
$$
where $$[b_2\bar\rho_1]= [b_2\bar\rho_1 g]= [\rho_2\bar{\hat
{b_2}}]=[\rho_2\bar{\hat {b_2}}g] =[5/2\rho].$$
\end{theorem}

\subsubsection{Nonlocal case: $k$ is even but not divisible by $4$, $k\neq 10$}

In this section we assume that $k$ is even but not divisible by $4$,
$k\neq 10$ and $\sigma$ is a dual GHJ subfactor appearing in $j\rho$
with $[\rho\bar\rho]=[0]+[k/2].$ By our assumption $H_{\sigma_1}$
must be of type $A$ or $D.$ If $H_{\sigma_1}$ is type $A$, then
multiply $\sigma_1$ on the right  by an automorphism if necessary,
we may assume that $\sigma_1\in \Delta,$ hence $\sigma_1\in
\Delta\rho.$ By Lemma \ref{end} either $1/2\sigma_1$ is irreducible
or $a_{1/2}\sigma_2$ is irreducible. Then  it is clear that the
possible pairs of $[\sigma_1,\sigma_2]$ are given by
$$
[i,\rho],i\neq [k/4],[1/2, b],[1/2, b'],
$$
where $[\rho a_{k/4}=[b]+b'], [1/2 b]=[1/2b']=[\rho
a_{k/4-1/2}].$\par If $H_{\sigma_1}$ is type $D$, we may assume that
$\sigma_1\in \Delta \rho,$ and it follows that $\sigma_2\in \bar\rho
\Delta\rho.$ It follows that the fusion graph of the action of
$a_{1/2}^{\rho_0}$ is given by one $D$ graph, whose vertices are
labeled by irreducible components of $a_{i}, i\in \Delta.$ By Lemma
\ref{end}, it follows that the following is a list of possible
intermediate subfactors: $ [\rho,w] $ where $w$ is an irreducible
component of $a_{i},i\neq k/4, [b, a_{1/2}], [ b,\tilde a_{1/2}]. $
\begin{theorem}\label{kodd}
Assume that  $k$ is even but not divisible by $4$, $k\neq 10$  and
$\sigma$ is a dual GHJ subfactor appearing in $j\rho$ with
$[\rho\bar\rho]=[0]+[k/2].$ Then the following is a list of possible
intermediate subfactors:
$$
[i,\rho],i\neq [k/4],[1/2, b],[1/2, b'],
$$
where $[\rho a_{k/4}]=[b]+b'], [1/2 b]=[1/2 b']=[\rho a_{k/4-1/2}],$
 $ [\rho,w] $ where $w$ is an irreducible component of
$a_{i},i\in \Delta, [ b, a_{1/2}], [ b,\tilde a_{1/2}]. $
\end{theorem}
\subsubsection {$k=10$}
By Lemma \ref{specd}, $H_{\sigma_1}$ may be $E_6$ when $\sigma=
3/2\rho_0.$ Suppose that $H_{\sigma_1}$ is $E_6.$ Then we can assume
that $\sigma_1\in \Delta \rho_2$ with $[\rho_2\bar\rho_2]=[0]+[3].$
It follows that $\sigma_2\in \bar\rho_2 \Delta \rho_0.$ But the
conjugate of $\bar\rho_2 \Delta \rho_0$ has already been determined
by Figure \ref{e6d}. By computing index we conclude that
$d_{\sigma_1}= d_{\rho_2}=d_{\sigma_2}.$  By equation (\ref{3/2}) we
conclude that $[\sigma_1,\sigma_2]=[\rho_2,\bar x] $ or
$[\sigma_1,\sigma_2]=[\rho_2,\bar{x} a_5].$

\begin{theorem}\label{10}
Assume that $k=10$ and $\sigma$ is a dual GHJ subfactor appearing in
$j\rho$ with $[\rho\bar\rho]=[0]+[5].$ Then the following is a list
of possible intermediate subfactors:
$$
[i,\rho],i\neq [k/4],[1/2,\rho b],[1/2,\rho b'],
$$
where $[\rho a_{k/4}]=[b]+[b'], [1/2\rho b]=[1/2\rho b']=[\rho
a_{k/4-1/2}],$
 $ [\rho,w] $ where $w$ is an irreducible component of
$a_{i},i\in \Delta,$ $ [ b, a_{1/2}], [ b,\tilde a_{1/2}] $ $
[\rho_2 , \bar{x}], [\rho_2 ,\bar {x} a_5]$ where
$$[\rho_2 \bar x]= [ \rho_2 \bar {x} a_5^{\rho}]
=[3/2\rho].$$
\end{theorem}

\subsection{$E_8$ case}

Assume that $\sigma$ is a dual GHJ subfactor appearing in $j\rho$
with $[\rho\bar\rho]=[0]+[5]+[9]+[14].$ Since $[8\rho]=[\rho],$ by
Lemma \ref{autosub} we may assume that $\rho=\rho_0\rho_1.$ Let
$\sigma=\sigma_1\sigma_2.$ If $H_{\sigma_1}$ is type $A$, then
multiply $\sigma_1$ on the right by an automorphism if necessary, we
may assume that $\sigma_1\in \Delta,$ hence $\sigma_1\in
\Delta\rho.$ By Lemma \ref{end} either $1/2\sigma_1$ is irreducible
or $a_{1/2}\sigma_2$ is irreducible. Then from figure \ref{e8} and
fusion rules it is clear that the possible pairs of
$[\sigma_1,\sigma_2]$ are given by
$$
[1/2,\rho],[1,\rho],[3/2,\rho],[2,\rho],[1/2,\rho b_3'], [1/2, b_4],
[1,\rho b_4]
$$
If  $H_{\sigma_1}$ is type $D$, we may assume that $\sigma_1\in
\Delta \rho_0,$ and it follows that $\sigma_2\in \bar\rho_0
\Delta\rho.$ From $\rho=\rho_0\rho_1$ we get $[\bar\rho_0\rho]
=[\rho_1]+[g\rho_1]$ where $[\bar\rho_0\rho_0]=[1]+[g], [g^2]=[0].$
It is then easy to determine   all irreducible sectors of
$\bar\rho_0 \Delta\rho.$ There are two $E_8$ graphs encoding the
action of $a_{1/2}^{\rho_0}$ in figure \ref{e8d}.
\begin{figure}[h]
\[
\xy (0,0)*{\bullet}; (10,0)*{\bullet}**\dir{-};
(20,0)*{\bullet}**\dir{-}; (30,0)*{\bullet}**\dir{-};
(40,0)*{\bullet}**\dir{-};(50,0)*{\bullet}**\dir{-};(60,0)*{\bullet}**\dir{-};
(40,0)*{}; (40,10)*{\bullet}**\dir{-}; (0,-3)*{\rho_1};
(10,-3)*{\rho_1 a_{1/2}}; (20,-3)*{\rho_1 a_{1}};(30,-3)*{\rho_1
a_{3/2}}; (40,-3)*{ \rho_1 a_2}; (50,-3)*{ \rho_1 b_3};
(40,13)*{\rho_1 b_3'};(60,-3)*{ \rho_1 b_4}
\endxy
\]

\[
\xy (0,0)*{\bullet}; (20,0)*{\bullet}**\dir{-};
(40,0)*{\bullet}**\dir{-}; (60,0)*{\bullet}**\dir{-};
(80,0)*{\bullet}**\dir{-};(100,0)*{\bullet}**\dir{-};(120,0)*{\bullet}**\dir{-};
(80,0)*{}; (80,10)*{\bullet}**\dir{-}; (0,-3)*{g\rho_1};
(20,-3)*{g\rho_1 a_{1/2}}; (40,-3)*{g\rho_1 a_{1}};(60,-3)*{g\rho_1
a_{3/2}}; (80,-3)*{ g\rho_1 a_2}; (100,-3)*{ g\rho_1 b_3};
(80,13)*{g\rho_1 b_3'};(120,-3)*{ g\rho_1 b_4}
\endxy
\]
\caption{Fusion graph of $a_{1/2}^{\rho_0}$}\label{e8d}
\end{figure}
By Lemma \ref{end} and figure \ref{e8d}, it follows that the
following is a list of possible intermediate subfactors: $
[\rho_0,w], $ where $w$ is one of the vertices in the first graph of
figure \ref{e8d},
$$
[\rho_0 b, \rho_1], [\rho_0 b, g\rho_1]
$$
where $\rho_0 b\prec [7\rho_0]$ and $[\rho_0 b\rho_1]= [\rho_0 b
g\rho_1]=[2\rho],$
$$
[1/2\rho_0,\rho_1],
[1\rho_0,\rho_1],[3/2\rho_0,\rho_1],[2\rho_0,\rho_1],[1/2\rho_0,\rho_1
b_3'], [1/2\rho_0,\rho_1 b_4], [1\rho_0,\rho_1 b_4]
$$

When $H_{\sigma_1}$ is $E_8,$ then there is $\rho_2\in H_{\sigma_1}$
such that $[\rho_2\bar\rho_2]=[\rho\bar\rho].$ By the cohomology
vanishing result of remark 5.4 in \cite{KLV}, multiplying $\sigma_1$
on the right by an automorphism if necessary, we can assume that
$[\rho_2]=[\rho],$ and so $\sigma_1\in \Delta \rho,$ and
$\sigma_2\in \bar \rho\Delta \rho.$ The set of irreducible sectors
of $\bar \rho\Delta \rho$ and fusion graphs of the action of
$a_{1/2}, \tilde a_{1/2}$ are given by  Figure 5 of \cite{BE3}. By
using Lemma \ref{end} and comparing indices it is straightforward to
determine the following list of possible intermediate subfactors:
$$
[\rho, a_{1/2}], [\rho, \tilde a_{1/2}], [\rho, a_{3/2}], [\rho,
\tilde a_{3/2}], [\rho, a_{1}], [\rho, \tilde a_{1}],
 [\rho, a_{2}], [\rho, \tilde a_{2}], [\rho, a_{1/2}\tilde{b_3'}],
$$
$$
[\rho, \tilde a_{1/2}b_3'],
 [\rho b_3', b_4], [\rho b_3',  a_{1/2}],[\rho b_3', \tilde a_{1/2}],
[\rho b_4, a_{1/2}], [\rho b_4,  \tilde a_{1/2}],[\rho b_4, a_{1}],
$$
$$
[\rho b_4,  \tilde a_{1}], [\rho b_4, b_3'], [\rho b_4,  \tilde
b_3'], [\rho a_{1/2}, b_4],  [\rho a_{1}, b_4],[\rho a_{1/2}, b_3'],
[\rho a_{1/2}, \tilde b_3']
$$
The additional fusion rules which are not immediate visible from
Figure \ref{e8} are
$$[b_3'b_4]=[a_{3/2}], [a_1b_4]=[a_2].$$

\begin{theorem}\label{casee8}
Assume that $\sigma$ is a dual GHJ subfactor appearing in $j\rho$
with $[\rho\bar\rho]=[0]+[5]+[9]+[14].$ Then the following is the
list of all intermediate subfactors that can occur:
$$
[1/2,\rho],[1,\rho],[3/2,\rho],[2,\rho],[1/2,\rho b_3'], [1/2, b_4],
[1,\rho b_4]
$$
$ [\rho_0,w], $ where $w$ is one of the vertices in the first graph
of figure \ref{e8d},
$$
[\rho_0 b, \rho_1], [\rho_0 b, g\rho_1]
$$
where $\rho_0 b\prec [7\rho_0]$ and $[\rho_0 b\rho_1]= [\rho_0 b
g\rho_1]=[2\rho],$
$$
[1/2\rho_0,\rho_1],
[1\rho_0,\rho_1],[3/2\rho_0,\rho_1],[2\rho_0,\rho_1],[1/2\rho_0,\rho_1
b_3'], [1/2\rho_0,\rho_1 b_4], [1\rho_0,\rho_1 b_4],
$$
$$
[\rho, a_{1/2}], [\rho, \tilde a_{1/2}], [\rho, a_{3/2}], [\rho,
\tilde a_{3/2}], [\rho, a_{1}], [\rho, \tilde a_{1}],
 [\rho, a_{2}], [\rho, \tilde a_{2}], [\rho, a_{1/2}\tilde{b_3'}],
 $$
 $$[\rho, \tilde a_{1/2}b_3'],
 [\rho b_3', b_4], [\rho b_3',  a_{1/2}],[\rho b_3', \tilde a_{1/2}],
[\rho b_4, a_{1/2}], [\rho b_4,  \tilde a_{1/2}],
$$
$$
[\rho b_4, a_{1}], [\rho b_4,  \tilde a_{1}], [\rho b_4, b_3'],
[\rho b_4, \tilde b_3'], [\rho a_{1/2}, b_4],  [\rho a_{1},
b_4],[\rho a_{1/2}, b_3'], [\rho a_{1/2}, \tilde b_3']
$$

The additional fusion rules which are not immediate visible from
Figure \ref{e8} are
$$[b_3'b_4]=[a_{3/2}], [a_1b_4]=[a_2].$$
\end{theorem}

\section{The lattice structure of intermediate dual GHJ subfactors}
In this section we list the lattices of intermediate subfactors of
dual GHJ subfactors. Given a dual GHJ subfactor $\sigma\prec \Delta
\rho,$ first we inspect all the pairs $[\sigma_1,\sigma_2]$ listed
in Th. \ref{casea}, Th. \ref{casee6}, Th. \ref{casee7},
Th.\ref{keven}, Th. \ref{kodd}, Th. \ref{16}, Th. \ref{10} and Th.
\ref{casee8} such that $[\sigma]=[\sigma_1\sigma_2].$ This gives all
intermediate subfactors of $\sigma.$ Then we use Cor. \ref{paircor}
and known fusion rules to determine the relations between these
intermediate subfactors. The result are listed in the following
figures. In each figure indexed by a dual GHJ subfactor $\sigma$ we
list all nontrivial intermediate subfactors $[\sigma_1,\sigma_2]$
with $[\sigma_1\sigma_2]=[\sigma]$. If $[\sigma_1,\sigma_2]$ lies
above $[\tau_1,\tau_2]$ and there is a line connecting them, then
$[\tau_1,\tau_2]\subset [\sigma_1,\sigma_2].$
\subsection{Type $A$}
When $k$ is odd or $k=2,$ all type $A$ GHJ subfactors are maximal.
When $k\geq 4$ is even, all $i\neq k/4$ are maximal.
\begin{figure}
\[\label{}
\xy (0,0)*{\bullet}; (0,-3)*{[\rho_0 b,\bar\rho_0]}
\endxy
\]
\caption{$1, k=4$ type $A$}\label{k4}
\end{figure}

\begin{figure}
\[\label{}
\xy (0,0)*{\bullet}; (0,10)*{\bullet}**\dir{-};
(0,-3)*{[\rho_0{b},\bar\rho_0]};(0,13)*{[\rho_0,b\bar\rho_0]}
\endxy
\]
\caption{$k/4, 4|k,k>4,$ type $A$}\label{a1}
\end{figure}

\par

\begin{figure}
\[
\xy (0,0)*{\bullet}; (20,0)*{\bullet};
(0,-3)*{[{b},\bar\rho_0]};(20,-3)*{[\rho_0,\bar b]}
\endxy
\]
\caption{$k/4$, $k\geq 6$ is even but not divisible by $4$, type
$A$}\label{a2}
\end{figure}

\subsection{Type $D$}
When $k=2$ the GHJ subfactor is maximal. When $k=4,$ $\rho_0 b,
\rho_0 b'$ are maximal. When $k$ is not divisible by $4,$ $b,b'\prec
k/4\rho_0$ are maximal. We note that $[\rho_0 b]=[\rho_0 b'g]$ when
$k$ is divisible by $4,$ so $\rho_0 b$ and $\rho_0 b'$ have
identical intermediate subfactor lattice.  The lattice in figure
\ref{dodd} for the case when $k=6$ was first obtained in \cite{GJ}
in the setting of type $II_1$ factors.

\begin{figure}
\[
\xy (0,0)*{\bullet};  (0,-3)*{[\rho_0,b]}
\endxy
\]
\caption{$\rho_0 b, 4|k, k\geq 8,$ type $D$}\label{d6}
\end{figure}

\begin{figure}
\[
\xy (0,0)*{\bullet}; (10,0)*{\bullet};
(0,-3)*{[i,\rho_0]};(10,-3)*{[\rho_0,a_i]}
\endxy
\]
\caption{$i\rho_0, i\neq k/4,k/4-1/2, 4|k, i\neq 5/2$ if $k=16,$
type $D$}\label{d1}
\end{figure}

\begin{figure}
\[
\xy (0,0)*{\bullet};
(0,20)*{\bullet}**\dir{-};(20,0)*{\bullet}**\dir{-};(20,20)*{\bullet};(40,0)*{\bullet};
(0,-3)*{[\rho,\bar{\hat{b_2}}]};(20,-3)*{[\rho,\bar{\hat{b_2}}g]};(0,23)*{[\rho_0,a_{5/2}]};
(0,20)*{\bullet}; (40,0)*{\bullet}**\dir{-}; (0,20)*{\bullet};
(60,0)*{\bullet}**\dir{-};
(20,23)*{[5/2,\rho_0]};(40,-3)*{[b_2,\bar\rho_1]};(60,-3)*{[b_2,\bar\rho_1g]}
\endxy
\]
\caption{$5/2\rho_0, k=16,$ type $D$}\label{d52}
\end{figure}

\begin{figure}
\[
\xy
(0,0)*{\bullet};(20,0)*{\bullet};(40,0)*{\bullet};(60,0)*{\bullet};
(0,-3)*{[1/2,\rho_0]}; (20,-3)*{[\rho_0 b,a_{1/2}]};
(40,-3)*{[\rho_0 b,\tilde a_{1/2}]}; (60,-3)*{[\rho_0,a_{1/2} ]}
\endxy
\]
\caption{$(1/2)\rho_0, [a_{1/2}b]=[a_{1/2}b']=[a_{1/2}],k=4,$ type
$D$}\label{d5}
\end{figure}

\begin{figure}[h]
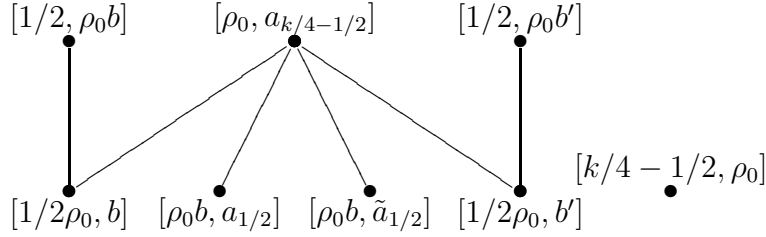

\[
\xy (0,0)*{\bullet}; (30,20)*{\bullet}**\dir{-};
(20,0)*{\bullet};(30,20)*{\bullet}**\dir{-};
(40,0)*{\bullet};(30,20)*{\bullet}**\dir{-};
(60,0)*{\bullet};(30,20)*{\bullet}**\dir{-};
(0,0)*{\bullet};(0,20)*{\bullet}**\dir{-};
(60,0)*{\bullet};(60,20)*{\bullet}**\dir{-}; (80,0)*{\bullet};
(0,-3)*{[1/2\rho_0,b]}; (20,-3)*{[\rho_0 b,a_{1/2}]};
(40,-3)*{[\rho_0 b,\tilde a_{1/2}]};(60,-3)*{[1/2\rho_0,b']};
(0,23)*{[1/2,\rho_0 b]}; (30,23)*{[\rho_0, a_{k/4-1/2}]};
(60,23)*{[1/2,\rho_0 b']}; (80,3)*{[k/4-1/2,\rho_0]}
\endxy
\]
\caption{$(k/4-1/2)\rho_0, [a_{1/2}b]=[a_{1/2}b']=[a_{k/4-1/2}],
4|k,k\geq 8,$ type $D$}\label{deven}
\end{figure}

\begin{figure}[h]
\[
\xy (0,0)*{\bullet}; (20,0)*{\bullet};
(0,-3)*{[i,\rho_0]};(20,-3)*{[\rho_0,a_i]}
\endxy
\]
\caption{$i\rho_0, i\neq k/4,k/4-1/2$, $k$ is even but not divisible
by $4, i\neq 3/2$ if $k=10,$ type $D$ }\label{d3}
\end{figure}

\begin{figure}[h]
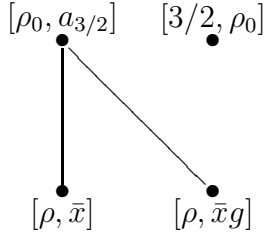

\[
\xy (0,0)*{\bullet};
(0,20)*{\bullet}**\dir{-};(20,0)*{\bullet}**\dir{-};(20,20)*{\bullet};
(0,-3)*{[\rho,\bar{x}]};(20,-3)*{[\rho,\bar{x}g]};(0,23)*{[\rho_0,a_{3/2}]};
(20,23)*{[3/2,\rho_0]}
\endxy
\]
\caption{$3/2\rho_0, k=10,$ type $D$}\label{d32}
\end{figure}

\begin{figure}
\[
\xy (0,0)*{\bullet}; (20,0)*{\bullet};(40,0)*{\bullet};
(60,0)*{\bullet};(80,0)*{\bullet}; (100,0)*{\bullet};
(0,-3)*{[\rho_0,a_{k/4-1/2}]}; (20,-3)*{[1/2, b]}; (40,-3)*{[1/2,
b']};(60,-3)*{[b,a_{1/2}]}; (80,-3)*{[b,\tilde
a_{1/2}]};(100,-3)*{[{k/4-1/2}, \rho_0]}
\endxy
\]
\caption{$(k/4-1/2)\rho_0, [a_{1/2}b]=[a_{1/2}b']=[a_{k/4-1/2}],$
$k\geq 6$ is not divisible by 4, }\label{dodd}
\end{figure}

\subsection{$E_6$}
In this case $\rho, \rho a_5$ are maximal. The lattice of
intermediate subfactors of $9/2\rho$ is isomorphic to that of
$1/2\rho$ since $[9/2\rho]=[1/2\rho a_5]$ and $a_5$ is an
automorphism.
\begin{figure}
\[
\xy (0,0)*{\bullet}; (20,0)*{\bullet};(40,0)*{\bullet};
(0,-3)*{[1/2,\rho]}; (20,-3)*{[\rho,a_{1/2}]};(40,-3)*{[\rho,\tilde
a_{1/2}]}
\endxy
\]
\caption{$1/2\rho,$ type $E_6$}\label{e61}
\end{figure}
\begin{figure}
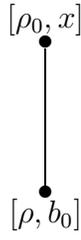

\[
\xy (0,0)*{\bullet}; (0,20)*{\bullet}**\dir{-}; (0,-3)*{[\rho,b_0]};
(0,23)*{[\rho_0,x]}
\endxy
\]
\caption{$\rho b_0,$ type $E_6$}\label{e62}
\end{figure}

\begin{figure}
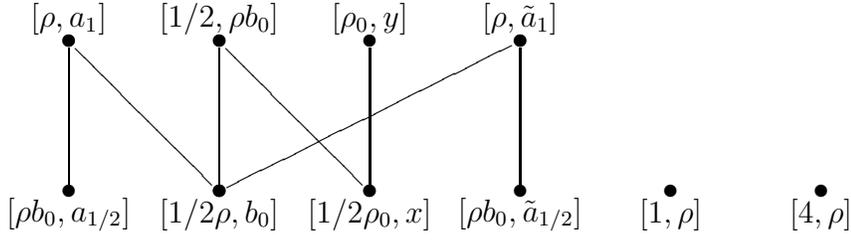

\[
\xy (0,0)*{\bullet}; (0,20)*{\bullet}**\dir{-};
(20,0)*{\bullet};(0,20)*{\bullet}**\dir{-};
(20,0)*{\bullet};(20,20)*{\bullet}**\dir{-};
 (20,0)*{\bullet};(60,20)*{\bullet}**\dir{-};
 (20,20)*{\bullet};(40,0)*{\bullet}**\dir{-};
 (40,0)*{\bullet};(40,20)*{\bullet}**\dir{-};
 (60,0)*{\bullet};(60,20)*{\bullet}**\dir{-};
 (80,0)*{\bullet}; (100,0)*{\bullet};
(0,-3)*{[\rho b_0,a_{1/2}]};(20,-3)*{[1/2\rho, b_0]}; (40,-3)*{[1/2
\rho_0 ,x]};(60,-3)*{[ \rho b_0 ,\tilde a_{1/2}]}; (0,23)*{[ \rho
,a_1]};(20,23)*{[ 1/2 ,\rho b_0]}; (40,23)*{[ \rho_0 ,y]};
(60,23)*{[ \rho ,\tilde a_1]};(80,-3)*{[1,\rho]};
(100,-3)*{[4,\rho]};
\endxy
\]
\caption{$1\rho ,$ type $E_6$}\label{e63}
\end{figure}



\begin{figure}
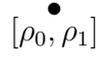

\[
\xy (0,0)*{\bullet}; (0,-3)*{[\rho_0,\rho_1]}
\endxy
\]
\caption{$\rho,$ type $E_7$}\label{e71}
\end{figure}

\begin{figure}
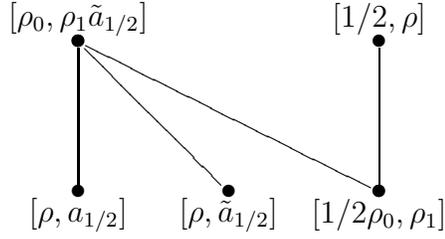

\[
\xy (0,0)*{\bullet};
(0,20)*{\bullet}**\dir{-};(0,20)*{\bullet};(20,0)*{\bullet}**\dir{-};
(0,20)*{\bullet}; (40,0)*{\bullet}**\dir{-};(40,0)*{\bullet};
(40,20)*{\bullet}**\dir{-}; (0,-3)*{[\rho ,a_{1/2}]};(20,-3)*{[\rho,
\tilde a_{1/2}]}; (40,-3)*{[1/2\rho_0 ,\rho_1]};(0,23)*{[ \rho_0
,\rho_1\tilde a_{1/2}]}; (40,23)*{[ 1/2 ,\rho]}
\endxy
\]
\caption{$1/2\rho ,$ type $E_7$}\label{e72}
\end{figure}

\begin{figure}
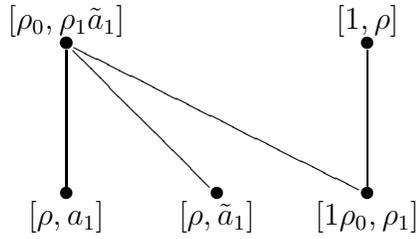

\[
\xy (0,0)*{\bullet};
(0,20)*{\bullet}**\dir{-};(20,0)*{\bullet}**\dir{-};
(0,20)*{\bullet}; (40,0)*{\bullet}**\dir{-};
(40,20)*{\bullet}**\dir{-}; (0,-3)*{[\rho ,a_{1}]};(20,-3)*{[\rho,
\tilde a_{1}]}; (40,-3)*{[1\rho_0 ,\rho_1]};(0,23)*{[ \rho_0
,\rho_1\tilde a_{1}]}; (40,23)*{[ 1 ,\rho]}
\endxy
\]
\caption{$1\rho ,$ type $E_7$}\label{e73}
\end{figure}

\begin{figure}
\[\xy (0,0)*{\bullet}; (0,-3)*{[\rho_0,\hat{b_1'}]}
\endxy
\]
\caption{$ b_1',$ type $E_7$}\label{e74}
\end{figure}

\clearpage

\begin{figure}[h]
\[
\xy (0,0)*{\bullet}; (0,-3)*{[\rho_0,\hat{b_2}]}
\endxy
\]
\caption{$ b_2,$ type $E_7$}\label{e75}
\end{figure}

\begin{figure}
\[\xy (0,0)*{\bullet};
(0,20)*{\bullet}**\dir{-};(20,0)*{\bullet}**\dir{-};(20,20)*{\bullet}**\dir{-};
(-20,0)*{\bullet}; (0,20)*{\bullet}**\dir{-};(-40,0)*{\bullet};
(0,20)*{\bullet}**\dir{-}; (0,-3)*{[\rho
,\tau]};(20,-3)*{[1/2\rho_0, \hat{b_2}]}; (-20,-3)*{[b_2,
,a_{1/2}]};(-40,-3)*{[ b_2 ,\tilde a_{1/2}]}; (0,23)*{[ \rho_0
,\hat{b_1}]};(20,23)*{[ 1/2 ,{b_2}]}
\endxy
\]
\caption{$b_1 ,$ type $E_7$}\label{e76}
\end{figure}

\begin{figure}
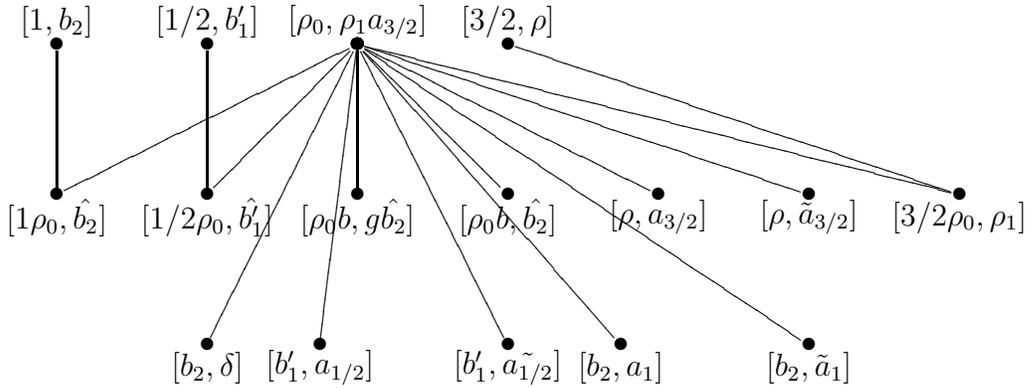

\[
\xy (0,20)*{\bullet};
(0,0)*{\bullet}**\dir{-};(40,20)*{\bullet}**\dir{-};
(20,20)*{\bullet};
(20,0)*{\bullet}**\dir{-};(40,20)*{\bullet}**\dir{-};
(40,0)*{\bullet};(40,20)*{\bullet}**\dir{-};
(60,0)*{\bullet};(40,20)*{\bullet}**\dir{-};
(80,0)*{\bullet};(40,20)*{\bullet}**\dir{-};
(100,0)*{\bullet};(40,20)*{\bullet}**\dir{-};(120,0)*{\bullet}**\dir{-};(60,20)*{\bullet}**\dir{-};
(20,-20)*{\bullet}; (40,20)*{\bullet}**\dir{-};
(35,-20)*{\bullet};(40,20)*{\bullet}**\dir{-}; (60,-20)*{\bullet};
(40,20)*{\bullet}**\dir{-};(75,-20)*{\bullet};
(40,20)*{\bullet}**\dir{-};(100,-20)*{\bullet};
(40,20)*{\bullet}**\dir{-}; (0,-3)*{[1\rho_0
,\hat{b_2}]};(20,-3)*{[1/2\rho_0, \hat{b_1'}]}; (40,-3)*{[\rho_0
b,g\hat{b_2}]};(60,-3)*{[ \rho_0 b ,\hat{b_2}]};(80,-3)*{[ \rho
,a_{3/2}]};(100,-3)*{[ \rho,\tilde a_{3/2}]};(120,-3)*{[ 3/2\rho_0
,\rho_1]}; (0,23)*{[ 1 ,{b_2}]};(20,23)*{[ 1/2 ,{b_1'}]}; (40,23)*{[
\rho_0,\rho_1 a_{3/2}]};(60,23)*{[ 3/2 ,\rho]}; (20,-23)*{[
b_2,\delta]}; (35,-23)*{[ b_1' ,a_{1/2}]};(60,-23)*{[ b_1'
,\tilde{a_{1/2}}]}; (75,-23)*{[ b_2 ,a_1]};(100,-23)*{[ b_2 ,\tilde
a_1]};
\endxy
\]
\caption{$3/2\rho ,$ type $E_7$}\label{e77}
\end{figure}


\begin{figure}
\[
\xy (0,0)*{\bullet}; (0,-3)*{[\rho_0,\rho_1]}
\endxy
\]
\caption{$\rho,$ type $E_8$}\label{e81}
\end{figure}

\begin{figure}
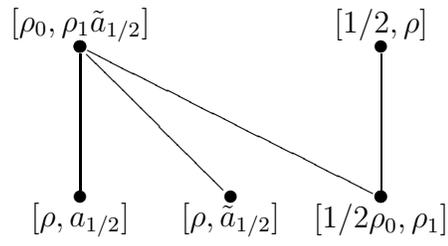

\[
\xy (0,0)*{\bullet};
(0,20)*{\bullet}**\dir{-};(0,20)*{\bullet};(20,0)*{\bullet}**\dir{-};
(0,20)*{\bullet}; (40,0)*{\bullet}**\dir{-};(40,0)*{\bullet};
(40,20)*{\bullet}**\dir{-}; (0,-3)*{[\rho ,a_{1/2}]};(20,-3)*{[\rho,
\tilde a_{1/2}]}; (40,-3)*{[1/2\rho_0 ,\rho_1]};(0,23)*{[ \rho_0
,\rho_1\tilde a_{1/2}]}; (40,23)*{[ 1/2 ,\rho]}
\endxy
\]
\caption{$1/2\rho ,$ type $E_8$}\label{e82}
\end{figure}

\begin{figure}
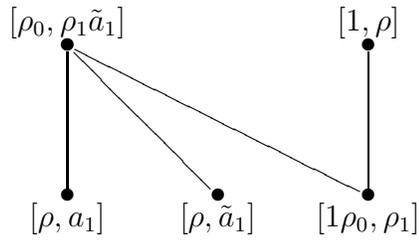

\[
\xy (0,0)*{\bullet};
(0,20)*{\bullet}**\dir{-};(20,0)*{\bullet}**\dir{-};
(0,20)*{\bullet}; (40,0)*{\bullet}**\dir{-};
(40,20)*{\bullet}**\dir{-}; (0,-3)*{[\rho ,a_{1}]};(20,-3)*{[\rho,
\tilde a_{1}]}; (40,-3)*{[1\rho_0 ,\rho_1]};(0,23)*{[ \rho_0
,\rho_1\tilde a_{1}]}; (40,23)*{[ 1 ,\rho]}
\endxy
\]
\caption{$1\rho ,$ type $E_8$}\label{e83}
\end{figure}

\begin{figure}
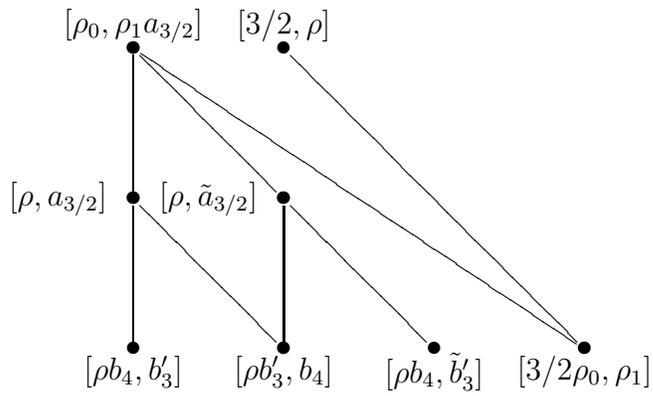

\[
\xy (0,0)*{\bullet}; (0,20)*{\bullet}**\dir{-}; (0,20)*{\bullet};
(0,40)*{\bullet}**\dir{-}; (20,0)*{\bullet};
(0,20)*{\bullet}**\dir{-};(20,0)*{\bullet};
(20,20)*{\bullet}**\dir{-};(20,20)*{\bullet};
(0,40)*{\bullet}**\dir{-};(20,20)*{\bullet};
(40,0)*{\bullet}**\dir{-};(0,40)*{\bullet};
(60,0)*{\bullet}**\dir{-};(60,0)*{\bullet};
(20,40)*{\bullet}**\dir{-}; (0,-3)*{[\rho b_4 ,b_3']};(20,-3)*{[\rho
b_3', b_4]}; (40,-3)*{[\rho b_4 ,\tilde b_3']};(60,-3)*{[ 3/2\rho_0
,\rho_1]}; (-10,20)*{[ \rho ,a_{3/2}]};(10,20)*{[\rho,\tilde
a_{3/2}] };(0,43)*{[\rho_0,\rho_1 a_{3/2}]};(20,43)*{[3/2,\rho]}
\endxy
\]
\caption{$3/2\rho ,$ type $E_8$}\label{e84}
\end{figure}

\begin{figure}
\[
\xy (0,0)*{\bullet};
(0,20)*{\bullet}**\dir{-};(40,0)*{\bullet}**\dir{-};(100,40)*{\bullet}**\dir{-};(40,20)*{\bullet}**\dir{-};
(20,40)*{\bullet}**\dir{-};(0,20)*{\bullet}**\dir{-};
(20,0)*{\bullet};(20,20)*{\bullet}**\dir{-};
(20,40)*{\bullet}**\dir{-};(20,20)*{\bullet};
(40,0)*{\bullet}**\dir{-}; (0,-3)*{[\rho b_4
,a_{1/2}]};(20,-3)*{[\rho b_4, \tilde a_{1/2}]}; (40,-3)*{[1/2\rho
,b_4]};(-10,20)*{[ \rho ,b_3]}; (10,20)*{[\rho,\tilde
b_3]};(55,20)*{[ 1/2\rho_0 ,\rho_1 b_4]};(20,43)*{[ \rho_0 ,\rho_1
b_3]};(100,43)*{[1/2  ,\rho b_4]}
\endxy
\]
\caption{$\rho b_3 ,$ type $E_8$}\label{e82n}
\end{figure}

\begin{figure}[h]
\[
\xy (0,0)*{\bullet}; (0,20)*{\bullet}**\dir{-};
(20,0)*{\bullet}**\dir{-}; (0,-3)*{[\rho
,b_3']};(20,-3)*{[\rho,\tilde{b_3'}]};(0,23)*{[\rho_0,\rho_1{b_3'}]}
\endxy
\]
\caption{$\rho b_3'$,type $E_8$}\label{e85}
\end{figure}

\begin{figure}[h]
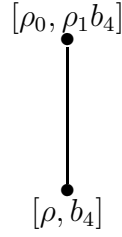

\[
\xy (0,0)*{\bullet}; (0,20)*{\bullet}**\dir{-};  (0,-3)*{[\rho
,b_4]};(0,23)*{[\rho_0,\rho_1{b_4}]}
\endxy
\]
\caption{$\rho b_4$,type $E_8$}\label{e86}
\end{figure}

\begin{figure}[h]
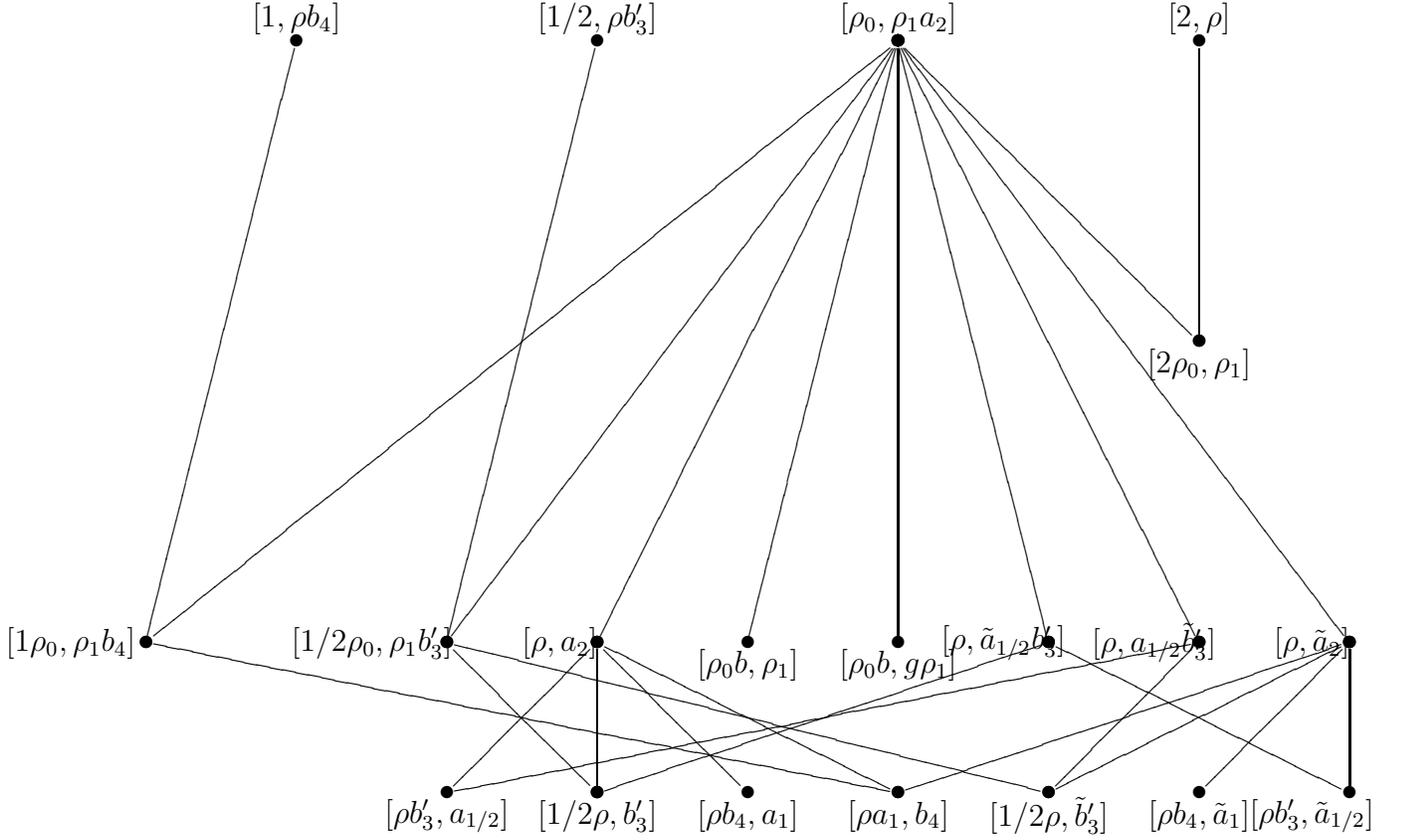

\[
\xy (0,80)*{\bullet}; (-20,0)*{\bullet}**\dir{-};
(80,80)*{\bullet}**\dir{-};(-20,0)*{\bullet};
(80,-20)*{\bullet}**\dir{-}; (40,80)*{\bullet};
(20,0)*{\bullet}**\dir{-}; (80,80)*{\bullet}**\dir{-};
(40,0)*{\bullet}**\dir{-};(40,-20)*{\bullet}**\dir{-};
(20,0)*{\bullet};(40,-20)*{\bullet}**\dir{-};
(40,0)*{\bullet};(60,-20)*{\bullet}**\dir{-};
(40,0)*{\bullet};(80,-20)*{\bullet}**\dir{-};(140,0)*{\bullet}**\dir{-};(120,-20)*{\bullet}**\dir{-};
(20,0)*{\bullet};(100,-20)*{\bullet}**\dir{-};
(120,0)*{\bullet};(100,-20)*{\bullet}**\dir{-};(140,0)*{\bullet}**\dir{-};
(60,0)*{\bullet};(80,80)*{\bullet}**\dir{-};
(80,0)*{\bullet};(80,80)*{\bullet}**\dir{-};
(100,0)*{\bullet};(80,80)*{\bullet}**\dir{-};
(120,0)*{\bullet};(80,80)*{\bullet}**\dir{-};
(140,0)*{\bullet};(80,80)*{\bullet}**\dir{-};
(120,40)*{\bullet};(80,80)*{\bullet}**\dir{-};
(100,0)*{\bullet};(40,-20)*{\bullet}**\dir{-};
(120,40)*{\bullet};(120,80)*{\bullet}**\dir{-}; (40,0)*{\bullet};
(20,-20)*{\bullet}**\dir{-}; (120,0)*{\bullet}**\dir{-};
(140,0)*{\bullet};
(140,-20)*{\bullet}**\dir{-};(100,0)*{\bullet}**\dir{-};
(20,-23)*{[\rho b_3',a_{1/2} ]};(135,-23)*{[\rho b_3',\tilde
a_{1/2}]}; (-30,0)*{[1\rho_0,\rho_1 b_4 ]};(10,0)*{[1/2\rho_0,\rho_1
b_3']}; (35,0)*{[\rho ,a_2]};(60,-3)*{[ \rho_0 b
,\rho_1]};(80,-3)*{[ \rho_0 b ,g\rho_1 ]};(94,0)*{[ \rho,\tilde
a_{1/2} b_3']};(114,0)*{[ \rho,a_{1/2} \tilde b_3']};(135,0)*{[
\rho,\tilde a_{2} ]}; (120,37)*{[ 2\rho_0,\rho_1]}; (40,-23)*{[
1/2\rho, b_3']};(60,-23)*{[ \rho b_4,a_{1}]}; (80,-23)*{[ \rho
a_{1},b_4]};(100,-23)*{[ 1/2\rho, \tilde b_3']};(120,-23)*{[ \rho
b_4,\tilde a_{1}]}; (0,83)*{[ 1,\rho b_4]};(40,83)*{[ 1/2,\rho
b_3']};(80,83)*{[ \rho_0,\rho_1 a_2]}; (120,83)*{[ 2,\rho ]};
\endxy
\]
\caption{$2\rho ,$ type $E_8$}\label{e88}
\end{figure}
\subsection{$E_7,E_8$ cases, and an example}
The labeling of type $E_7, E_8$ cases are described in Th.
\ref{casee7} and Th. \ref{casee8} respectively. Here we use the most
complicated case of $2\rho ,$ type $E_8$ to explain how we obtain
the lattice structure. Let us first use Lemma \ref{pair} to show
that $[\rho,a_2]\neq [\rho,\tilde a_2].$ If $[\rho,a_2]=
[\rho,\tilde a_2],$ we can find an automorphism $\sigma$ such that
$[\rho\sigma]=[\rho], [\sigma a_2]=[\tilde a_2].$ But $\tilde a_2
a_2$ is irreducible, and this implies that $[\sigma]=[\tilde a_2
a_2],$ contradicting our assumption that $\sigma$ is an
automorphism. Alternatively we can argue from $[\rho\sigma]=[\rho]$
that $\sigma\prec \bar\rho \rho,$ and from the formula for $\bar\rho
\rho$ in \cite{BEK1} we conclude that $[\sigma]=[0],$ and hence
$[a_2]=[\tilde a_2],$ again a contradiction. A good way to look at
the intermediate subfactors of $2\rho ,$ type $E_8$ is to start with
dual GHJ of type $A$: these are pairs with the first components
labeled by a half integer, and there are three of them; for type $D$
the first components are labeled by a half integer multiplied by
$\rho_0,$ and there are six of them; for type $E$ the first
components are labeled the vertices of figure \ref{e8},  and there
are eleven of them.\par Let us now explain why $[\rho
b_4,a_1]\subset [\rho,a_2]$ and yet $[\rho b_4,a_1]$ is not a
subfactor of $[\rho,\tilde a_2].$ Since by Lemma \ref{pair}, $[\rho
b_4,a_1]\subset [\rho,b_4 a_1],$ but $[b_4a_1]=[a_2],$ so we have
shown that  $[\rho b_4,a_1]\subset [\rho,a_2].$ Now if $[\rho
b_4,a_1]\subset [\rho,\tilde a_2],$ by Lemma \ref{pair} we can find
$\sigma$ such that
$$[\rho b_4]=[\rho \sigma], [\sigma a_1]=[\tilde a_2]$$
It follows that $[\rho,\sigma]$ is an intermediate subfactor of
$\rho b_4.$ But all such pairs are classified in Th. \ref{casee8},
and by inspection we conclude that $[\rho,\sigma]=[\rho,b_4].$ By
definition we have an automorphism $\sigma_1$ such that
$[\rho\sigma_1]=[\rho], [\sigma_1\sigma]=[b_4].$ Since by
\cite{BEK1} the only subsector of $\bar\rho\rho$ which is an
automorphism is $[0],$ we have $[\sigma]=[b_4],$ and $[\sigma
a_1]=[a_2]=[\tilde a_2],$ a contradiction. The rest of relations in
Figure \ref{e88} are derived in a similar way.

\subsection{A negative answer to question \ref{gag}}\label{neg}
From the list of lattices we can easily verify conjecture \ref{wall}
for GHJ subfactors. It is an interesting question to check that
whether the stronger conjecture \ref{mod} is true for all GHJ
subfactors using our list. \par We  claim that the  GHJ subfactor $
\bar \rho 3/2 $ of type $E_7$ gives a negative answer to question
\ref{gag} in the appendix. This subfactor is dual to $3/2\rho $
whose lattice of intermediate subfactors is given by figure
\ref{e77}. Since $3/2\rho$ has $12$ minimal subfactors, it follows
that $\bar\rho 3/2$ has $12$ maximal subfactors. Note that
$[\bar\rho 3/2][3/2\rho]\in \bar\rho\Delta\rho,$ and it is known by
Figure 42 of \cite{BEK2} that the only sector  with index $1$ which
appears in $\bar\rho\Delta\rho$ is $[0].$ By Lemma \ref{autosub1} we
have $\Auto(\bar\rho 3/2)$ is trivial. On the other hand it is easy
to calculate $[3/2\rho\bar\rho 3/2],$ and we find there are $9$
irreducible sectors which can appear in $[3/2\rho\bar\rho 3/2].$
Since $12>9,$ this gives a negative answer to problem \ref{gag}.

\appendix
\section{A proof of relative version of Wall's conjecture for
solvable groups} As pointed in introduction, Wall proved his
conjecture for solvable group.  Another supporting evidence, as
observed by V. F. R. Jones, is that the minimal version of
conjecture  \ref{wall} holds for subfactors with a commutative
$N'\cap M_1.$ In this appendix we will prove the relative version of
Wall's conjecture for solvable groups. Our proof is partially
inspired by some comments of V. F. R. Jones. Let $N\subset M$ be an
irreducible subfactor with finite Jones index, $e_N$ the Jones
projection from $M$ to $N$, and let $P_i, 1\leq i\leq n$ be the set
of minimal intermediate subfactors. Then the Jones projections $e_i$
from $M$ onto $P_i$ are in $N'\cap M_1.$ Conjecture \ref{wall} will
follow if $e_i,1\leq i\leq n, e_N$ are linearly independent.
Unfortunately this is not true in general. In fact, Let $G$ acts
properly on the hyperfinite type $II_1$ factor $R,$ and consider the
subfactors $R^G\subset R^H\subset R$ where $R^G, R^H$ are fixed
point subalgebras of $R$ under the action of $G, H$ respectively.
Let $K_i, 1\leq i\leq n$ be the set of maximal subgroups of $G$
which strictly contains $H.$ Let $e_i$ be the Jones projections from
$R$ onto $R^{K_i}.$ Note that $e_i =\frac{1}{|K_i|}\sum_{g\in K_i}
g\in {\Bbb C G}.$ For any subgroup $K$ of $G$ we will denote by
$e_K=\frac{1}{|K|}\sum_{g\in K} g\in {\Bbb C G}.$ Note that ${\Bbb C
G}$ is a $C^*$ algebra. We denote by $l(G)$ the abelian algebra of
complex valued functions on $G.$ There are examples when $G$  is a
semidirect product of an elementary abelian group $V$ by $G_1$ which
acts irreducibly on $V,$ and we find that the set
$\frac{1}{|G_1|}\sum_{g\in G_1} vgv^{-1}, v\in V$ is not linearly
independent (Note that for any $v\in V,$ $vG_1v^{-1}$ is maximal in
$G$ by our assumption). However the following modification seems to
be interesting:
\begin{conjecture}\label{mod}
Let $N\subset M$ be an irreducible subfactor with finite Jones
index, and let $P_i, 1\leq i\leq n$ be the set of minimal
intermediate subfactors. Denote by $e_i\in N'\cap M_1, 1\leq i\leq
n$ the Jones projections $e_i$ from $M$ onto $P_i$ and $e_N$ the
Jones projections $e_N$ from $M$ onto $N.$ Then there are vectors
$\xi_i, \xi\in N'\cap M_1$ such that $e_i \xi_i=\xi_i,1\leq i\leq n,
e_N\xi=\xi,$ and $\xi_i, 1\leq i\leq n,\xi$ are linearly
independent.
\end{conjecture}
\begin{remark}
We note that unlike conjecture \ref{wall}, the conjecture above
makes use of the algebra structure of $N'\cap M_1$ and therefore
does not immediately imply the dual version or if one replaces
minimal by maximal.
\end{remark}
By definition conjecture \ref{mod} implies conjecture \ref{wall},
and we shall prove conjecture \ref{mod} for $R^G\subset R, G$
solvable. In fact it is easy to check that for $R^G\subset R$
conjecture \ref{mod} is equivalent to:
\begin{conjecture}\label{mod2}
Let $K_i, 1\leq i\leq n$ be a set of maximal subgroups of $G.$ Then
there are vectors $\xi_i\in l(G), 1\leq i\leq n$ such that $
e_G\xi_i=0,$ $\xi_i$ are $K_i$ invariant and linearly independent.
\end{conjecture}
We will prove conjecture \ref{mod2} when $G$ is solvable.
\begin{lemma}\label{prim}
Suppose that $K$ acts irreducibly on an elementary abelian group
$V,$ and $G$ is the semi-direct product of $V$ by $K.$  For each $v
K v^{-1}, v\neq 0,$ we assign a vector
$\xi_v:=\delta_{vK}-\frac{1}{|V|}1\in l(G),$ and for $K,$ we assign
$\xi_K:=\delta_{Mb}-\frac{1}{|V|}1\in l(G)$ where $b\neq 0,$ and we
use $\delta_S$ to denote the characteristic function of a set
$S\subset G,$ and $1$ stands for constant function with value $1.$
Then $\xi_v, v\in V, v\neq 0, \xi_K$ verify conjecture \ref{mod2}
for $vKv^{-1}, v\in V.$
\end{lemma}
\prf By definition we just have to check that $\xi_v, v\in V, v\neq
0, \xi_K$ are linearly independent. Suppose that $\lambda_v$ are
complex numbers such that
$$
\sum_{v\neq 0,v\in V} \lambda_v \xi_v + \lambda_0 \xi_K =0.
$$
The we have
$$
\sum_{v\neq 0,v\in V} \lambda_v \delta_{vK} + \lambda_0 \delta_{Kb}-
\sum_{v\in V} \lambda_v \frac{1}{|V|} 1 =0.
$$
For a fixed $v\in V,v\neq 0,$ since $K$ acts irreducibly on $V,$ we
can find $k\in K$ such that $k(v):=k^{-1} v k\neq b.$ It follows
that $\delta_{Kb}(vk)=0.$  Evaluate the LHS of the above sum at $vk,
v\neq 0$ we have
$$
\lambda_v= \frac{1}{|V|}\sum_{v'\in V} \lambda_{v'}.
$$
Since $k b= k^{-1}(b) k,$ evaluating the above sum at $kb$ we have
$$
\lambda_0 + \lambda_{k^{-1}(b)}-\frac{1}{|V|}\sum_{v\in V}
\lambda_v=0.
$$
Notice that $k^{-1}(b)\neq 0$ since $b\neq 0,$ we conclude that
$\lambda_0=0,$  all $\lambda_v, v\neq 0$ are identical to the same
value $\lambda,$ and $\frac{1}{|V|}(|V|-1)\lambda=\lambda.$ It
follows that $\lambda_v=0, \forall v\in V.$ \qed
\begin{lemma}\label{induction}
Suppose that $H$ is a normal subgroup of $G$ and $K_i\geq H, 1\leq
i\leq n$ is a set of maximal subgroups of $G.$ If conjecture
\ref{mod2} is true for $K_i/H\leq G/H,$ then it is also true for
$K_i\leq G, 1\leq i\leq n.$
\end{lemma}
\prf By assumption we have functions $\xi_i: G/H\rightarrow \C$ such
that $\xi_i $ is linearly independent and $K_i/H$ invariant,
$e_{G/H} \xi_i=0, 1\leq i\leq n.$  Let $\pi: G\rightarrow G/H$ be
the projection, then $\xi_i\pi: G\rightarrow \C$ is linearly
independent and $K_i$ invariant,$e_{G} \xi_i\pi=0, 1\leq i\leq n.$
\qed

Recall that for a subgroup $K\leq G,$ $\Core_G K:=\cap_{g\in G}
gKg^{-1}$ is the largest normal subgroup of $G$ that is contained in
$K.$

\begin{lemma}\label{es}
Let $H_1,H_2$ be subgroups of a finite group $G.$ Denote by $H_1H_2$
the set of different elements $g$ which can be written as $h_1h_2$
with $h_1\in H_1, h_2\in H_2.$ Then
$e_{H_1}e_{H_2}=\frac{1}{|H_1H_2|}\sum_{g\in H_1H_2} g.$
\end{lemma}
\prf Let $g=h_1h_2$ with $h_1\in H_1, h_2\in H_2.$ Then
$h_1h_2=h_1'h_2'$ iff $h_1^{-1}h_1'=h_2{h_2'}^{-1}\in H_1\cap H_2.$
It follows that for each $g=h_1h_2\in H_1H_2,$ there are $H_1\cap
H_2$ different pairs of $(h_1',h_2')\in H_1\times H_2$ such that
$g=h_1'h_2'.$ Hence $|H_1H_2|=\frac{|H_1||H_2|}{|H_1\cap H_2|}$ and
the lemma follows from definition. \qed

\begin{proposition}\label{sol}
Conjecture \ref{mod2} is true for $G$ solvable.
\end{proposition}
\prf The proof goes by induction on $|G|n.$ Consider $H:=\Core_G
K_1.$ If $\Core_G K_i$ does not contain $H,$ then $K_i H\neq K_i,$
and since $K_i H$ is a subgroup of $G, K_i$ is maximal, we have $K_i
H=G.$ Suppose that there is at least one $K_i$ such that  $\Core_G
K_i$ does not contain $H.$ By  induction hypothesis, for the set of
$K_i$ with $\Core_G K_i$ not containing $H,$ we can find vectors
$\xi_i$ which verifies conjecture \ref{mod2}, and  for the set of
$K_i$ with $\Core_G K_i$  containing $H,$ we can find vectors
$\xi_i$ which verifies conjecture \ref{mod2}. We claim that such set
$\xi_i, 1\leq i\leq n $  is linearly independent.  Suppose that
$\lambda_i\in \C,1\leq i\leq n $ such that $\sum_{1\leq i\leq n }
\lambda_i \xi_i=0.$ Multiply on the left by $e_{H}.$ We have
$$
\sum_{i, \Core_G K_i\geq H} \lambda_i e_H \xi_i +\sum_{i, \Core_G
K_i\cap H\neq H} \lambda_i e_H \xi_i=0.
$$
Note that if $ \Core_G K_i\geq H,$ then $e_H \xi_i= e_H
e_{K_i}\xi_i=e_{K_i}\xi_i;$ if $\Core_G K_i$ does not contain $H,$
then $e_H \xi_i= e_H e_{K_i}\xi_i=e_G\xi_i=0$ by Lemma \ref{es}. It
follows that
$$
\sum_{i, \Core_G K_i\geq H} \lambda_i  \xi_i=0, \sum_{i,
\Core_GK_i\cap H\neq H} \lambda_i  \xi_i=0.
$$
By our assumption $\lambda_i=0,1\leq i\leq n. $\par

So we are left with the case that  $\Core_G K_i \geq H,1\leq i\leq
n. $ By replacing $H=\Core_G K_1$ by $H=\Core_G K_j $ for some
$1\leq j\leq n$ we can now assume that $\Core_G K_i=H,1\leq i\leq n.
$ If $H$ is nontrivial, by Lemma \ref{induction} we are done.\par

If $H$ is trivial, by Th. 15.6 of \cite{solvable}, $G$ is the
semidirect product of an elementary abelian group $V$ by $K_1,$ and
the action of $K_1$ on $V$ is irreducible. Moreover by Th.16.1 of
\cite{solvable} all maximal subgroup $K$ of $G$ with $\Core_G K=H$
is of the form $vK_1v^{-1}$ for some $v\in V.$ By Lemma \ref{prim}
we are done. \qed

\begin{remark}\label{coh}
The reduction in the proof of the above proposition works for
general groups, and conjecture \ref{mod2} can be reduced to the case
where $G$ is a primitive group, and the set of maximal subgroups
have trivial core. Such groups are classified by O'Nan-Scott theorem
(cf. \S 4 of \cite{NS}). The first case is when  $G$ is the
semidirect product of an elementary abelian group $V$ by $K_1,$ and
the action of $K_1$ on $V$ is irreducible.  When $G$ is not
solvable,  maximal subgroups $K$ of $G$ with trivial core are not
conjugates of $K_1$, and our proof above does not work. Such maximal
subgroups are related to the first cohomology of $K_1$ with
coefficients in $V,$ and conjecture \ref{mod2} implies that the
order of this cohomology is less than $|K_1|$ (cf. Question 12.2 of
\cite{lub}). Unfortunately even though it is believed that the order
of this cohomology is small (cf. \cite{Gu}), the bound $|K_1|$ has
not been achieved yet.
\end{remark}
\begin{corollary}\label{relwall}
The relative version of Wall's conjecture is true for solvable
groups.
\end{corollary}
\prf Let $K_i, 1\leq i\leq n$ be a set of maximal subgroups of $G$
strictly containing $H.$ By Prop. \ref{sol} we can find vectors
$\xi_i\in l(G), 1\leq i\leq n$ such that $ e_G\xi_i=0,$ $\xi_i$ are
$K_i$ invariant and linearly independent. Since $K_i\geq H,$ we have
$e_H\xi_i= e_H e_{K_i}\xi_i= e_{K_i}\xi_i=\xi_i, 1\leq i\leq n.$ It
follows that $\xi_i$ is $H$ invariant, and can be thought as
functions on $l(G/H).$ Since $e_G\xi_i=0,$ we conclude that $1,
\xi_i, 1\leq i\leq n$ are linearly independent functions on
 $l(G/H)$ and the corollary follows.
\qed

At the end of this appendix we discuss a question which is motivated
by the following conjecture of Aschbacher-Guralnick in \cite{AsbG}:
\begin{conjecture}\label{ag}
Let $G$ be a finite group. Then the number of conjugacy classes of
maximal subgroups is less or equal to the number of conjugacy
classes of $G.$
\end{conjecture}
Conjecture \ref{ag} was proved in \cite{AsbG} for solvable $G.$ Here
we give a slightly proof (with strict inequality) in the spirit of
proof of Prop. \ref{sol}.
\begin{proposition}\label{strict}
If $G$ is a finite solvable group, then  the number of conjugacy
classes of maximal subgroups is less than the number of conjugacy
classes of $G.$
\end{proposition}

\prf Let $K_i, 1\leq i\leq n$ be a set of representatives of
conjugacy classes of maximal subgroups of $G.$ Let $C_1,...,C_k$ be
the conjugacy classes of $G.$ Define $f_i:=\sum_{g\in G}
ge_{K_i}g^{-1} - |G| e_G, 1\leq i\leq n, $ $h_j:=\sum_{g\in C_j} g,
1\leq j\leq k.$ Note that $e_Gf_i=0, e_G=|G|^{-1}\sum_{1\leq j\leq
k} h_j, $ $h_j, 1\leq j\leq k$ is linearly independent, and $f_i$ is
in the space spanned by $h_j, 1\leq j\leq k.$  We claim that $f_i,
1\leq i\leq n$ are linearly independent. This will prove the
proposition since
 $f_i$ is in the space
spanned by $h_j, 1\leq j\leq k,$ and $\sum_{1\leq j\leq k} h_j
f_i=0.$\par

Suppose that $\lambda_i\in \C,1\leq i\leq n $ such that $\sum_{1\leq
i\leq n } \lambda_i f_i=0.$ Fix $g$ and $1\leq i\leq n.$ If $j\neq
i,$ then $gK_ig^{-1}$ is not conjugate to $hK_jh^{-1},h\in G$ and by
Th. 16.2 of \cite{solvable} we have $gK_ig^{-1}hK_jh^{-1}=G,$ and by
Lemma \ref{es} $ge_{K_i}g^{-1} he_{K_j}h^{-1}= e_G.$  Multiply
$\sum_{1\leq i\leq n } \lambda_i f_i=0$ on the left by
$ge_{K_i}g^{-1}$ we have
$$
\lambda_i ge_{K_i}g^{-1} f_i= ge_{K_i}g^{-1}\lambda_i(\sum_{h\in G}
he_{K_i}h^{-1} - |G| e_G) =0, \forall g\in G.
$$
It follows that $|\lambda_i|^2 f_i f_i^*=0,$ and since $\C G$ is a
$C^*$ algebra,  $\lambda_i f_i=0.$ By looking at the coefficient of
identity element of $G$ in $f_i$ we conclude that $\lambda_i
(|G|/|K_i|-1)=0.$ since $|G|> |K_i|$ we conclude that $\lambda_i=0.$
\qed

The following question is motivated by  conjecture \ref{ag}:
\begin{question}\label{gag}
Let $N\subset M$ be an irreducible subfactor with finite index. Let
$\Auto (M|N):=\{ \alpha\in \Auto(M)| \alpha(n)=n, \forall n\in N\}.$
We say that two intermediate subfactors $P_1,P_2$ are conjugate if
there is an $\alpha\in \Auto(M|N)$ such that $\alpha(P_1)=(P_2).$ Is
the number of conjugacy classes of maximal or minimal subfactors
less or equal to the number of irreducible representations of
$N'\cap M_1$ ?
\end{question}
\begin{remark}
There is a similar formulation of the above question using $\Normal
(M|N)$ and Lemma \ref{normalsub1}.
\end{remark}
Take $N=M^G\subset M.$ Then $\Auto(M|M^G)=G, N'\cap M_1=\C G.$ The
conjugacy classes of maximal  subfactors is the same as the
conjugacy classes of minimal subgroups of $G,$ and it is easy to see
that the number conjugacy classes of  minimal subgroups of $G$ is
less than the number of conjugacy class of $G$, which is the same as
the number of irreducible representations of $N'\cap M_1=\C G.$ On
the other hand the conjugacy classes of minimal  subfactors is the
same as the conjugacy classes of maximal subgroups of $G,$ and
question  \ref{gag} is equivalent to conjecture \ref{ag}.\par

Let $M$ be the cross product of $N$ by $G.$ Then
 $\Auto(M|N)$ is isomorphic to the set of one dimensional representations of $G,$ and
 $N'\cap M_1=l(G).$ In this case $\Auto(M|N)$ preserves every
intermediate subfactor.

The conjugacy classes of minimal  subfactors can be identified  as
the set of minimal subgroups of $G,$ and it is easy to see that the
number of  minimal subgroups of $G$ is less than $|G|,$ which is the
same as the number of irreducible representations of $N'\cap
M_1=l(G).$ On the other hand the conjugacy classes of maximal
subfactors can be identified  as the set of  maximal subgroups of
$G,$ and question \ref{gag} is equivalent to Wall's conjecture (with
$\leq$ instead of $<$).

However, as shown in \S\ref{neg} question \ref{gag} has a negative
answer for general subfactors. Is there any natural modification of
the statement in question \ref{gag} so that it has a better chance
of being true while still generalizing conjecture \ref{ag}?

\end{document}